\documentclass[11pt,twoside]{amsart}

\usepackage{amsmath}
\usepackage{amsfonts}
\usepackage{amssymb}

\usepackage[all]{xy}

\title
{Symmetry of the definition  of degeneration in triangulated categories}

\author{Manuel Saor\'\i{}n}
\address{\newline
Departemento de Matem\'aticas,
\newline Universidad de Murcia, Aptdo. 4021
\newline 30100 Espinardo, Murcia,
\newline Spain}
\email{msaorinc@um.es}

\author{Alexander Zimmermann}
\address{\newline
Universit\'e de Picardie,
\newline D\'epartement de Math\'ematiques et LAMFA (UMR 7352 du CNRS),
\newline 33 rue St Leu,
\newline F-80039 Amiens Cedex 1,
\newline France}
\email{alexander.zimmermann@u-picardie.fr}

\newtheorem{Lemma1}{{Lemma}}
\newtheorem{Theo1}[Lemma1]{{Theorem}}
\newtheorem{Def1}[Lemma1]{{Definition}}
\newtheorem{Prop1}[Lemma1]{{Proposition}}
\newtheorem{Claim1}[Lemma1]{{Claim}}
\newtheorem{Rem1}[Lemma1]{{Remark}}
\newtheorem{Cor1}[Lemma1]{{Corollary}}
\newtheorem{Ex1}[Lemma1]{{Example}}
\newtheorem{Qu1}[Lemma1]{{Question}}

\newenvironment{Lemma}{\begin{Lemma1}}{\end{Lemma1}}
\newenvironment{Def}{\begin{Def1}\em}{\end{Def1}}
\newenvironment{Prop}{\begin{Prop1}}{\end{Prop1}}
\newenvironment{Rem}{\begin{Rem1}\rm}{\end{Rem1}}
\newenvironment{Theorem}{\begin{Theo1}}{\end{Theo1}}

\newenvironment{Cor}{\begin{Cor1}}{\end{Cor1}}

\newenvironment{Example}{\begin{Ex1}\em}{\end{Ex1}}

\newcommand{\uar}{\uparrow}

\newcommand{\lra}{\longrightarrow}

\newcommand{\ra}{\rightarrow}
\newcommand{\sdp}{\times\kern-.2em\vrule height1.1ex depth-.05ex}
\newcommand{\epi}{\lra \kern-.8em\ra}

\newcommand{\N}{{\mathbb N}}

\newcommand{\B}{{\mathbb B}}

\newcommand{\ol}{\overline}

\newcommand{\Z}{{\mathbb Z}}

\newcommand{\cone}{\text{cone}}

 \normalbaselines


\thanks{The first named author is backed by research projects from the Ministerio de Econom\'ia y Competitividad of Spain (MTM201346837-P) and the Fundaci\'on `S\'eneca' of Murcia (19880/GERM/15), both with a part of FEDER funds. He thanks these institutions for their support.}

\date{December 16, 2016}

\begin{document}

\begin{abstract}
Module structures of an algebra on a fixed finite dimensional vector space form
an algebraic variety. Isomorphism classes correspond to orbits of the action of
an algebraic group on this variety and a module is a degeneration of another if it belongs
to the Zariski closure of the orbit. Riedtmann and Zwara gave an algebraic
characterisation of this concept in terms of the existence of short exact
sequences. Jensen, Su and Zimmermann, as well as independently Yoshino,  studied
the natural generalisation of the Riedtmann-Zwara degeneration to triangulated
categories. The definition has an intrinsic non-symmetry. Suppose that we have a
triangulated category in which idempotents split and either for which the endomorphism
rings of all objects are artinian, or which is the category of compact objects in an
algebraic compactly generated triangulated $K$-category. Then we show that the non-symmetry in the
algebraic definition of the degeneration is inessential in the sense that the
two possible choices which can be made in the definition lead to the same concept.
\end{abstract}

\maketitle

\section*{Introduction}

For a finite dimensional $K$-algebra over an algebraically closed field
$K$ the set of $d$-dimensional $A$-modules is just the space
of $K$-algebra homomorphisms from $A$ to the algebra of $d$ by $d$
matrices over $K$. It
carries therefore the structure of an algebraic variety
$mod(A,d)$, and allows a
$GL_d(K)$-action given by conjugation of matrices.
$GL_d(K)$-orbits correspond to isomorphism
classes of modules, and we say that a $d$-dimensional
module $M$ corresponding to the point
$m\in mod(A,d)$ degenerates to the module $N$ with
corresponding point $n\in mod(A,d)$ if $n$
belongs to the Zariski-closure of the orbit $GL_d(K)\cdot m$.
We write in this case $M\leq_{\text{deg}}N$. It is clear that $\leq_{\text{deg}}$ is a partial
order on the set of isomorphism classes of finite dimensional $A$-modules.
Zwara and Riedtmann defined another relation between $A$-modules, namely
$M\leq_{\text{Zwara}} N$ if and only if there is a finite
dimensional $A$-module $Z$ and a short exact sequence $0\ra N\ra M\oplus Z\ra Z\ra 0$.
Moreover, they showed in
\cite{Riedtmann1,Zwara1}
$$M\leq_{\text{deg}}N\Leftrightarrow M\leq_{\text{Zwara}} N.$$
Zwara showed in \cite{Zwara} by a purely algebraic arguments that in the
category of finite dimensional modules over an algebra
there is $Z$ and a short exact sequence $0\ra N\ra M\oplus Z\ra Z\ra 0$
if and only if there is $Z'$ and a short exact sequence
$0\ra Z'\ra M\oplus Z'\ra N\ra 0$.

In joint work with Jensen and Su \cite{JSZ}, and independently
by Yoshino in \cite{Yoshino} for the (triangulated)
stable category of maximal Cohen Macaulay modules
over local Gorenstein $k$-algebras,
the concept $\leq_{\text{Zwara}}$ was generalised
in the obvious way to general triangulated categories. More precisely, for a
triangulated category $\mathcal T$ we define for two objects $M$ and $N$ that
$M\leq_\Delta N$ if and only if there is an object $Z$ and a distinguished triangle
$Z\stackrel{u\choose v}{\lra} M\oplus Z\stackrel{}{\lra} N\ra Z[1]$.
Yoshino insisted in the point that one should ask that the
induced endomorphism $v$ of $Z$ is nilpotent. Using Fitting's lemma and possibly replacing
$Z$ by a suitable direct summand,  this is automatic if one
assumes Krull-Schmidt properties and artinian endomorphism rings for all objects.
We denote $M\leq_{\Delta+\text{nil}}N$ if
$M\leq_\Delta N$ and the induced endomorphism $v$ on $Z$ is nilpotent.
The concept $\leq_\Delta$ was used in an essential way in work of Keller
and Scherotzke on Nakajima quiver varieties.
\cite{JSZ} concentrated on partial order properties of $\leq_\Delta$.
Further conditions guaranteeing partial order properties of $\leq_\Delta$ can
be found for various situations in \cite{Webb}, \cite{Wang} and \cite{SZ}.
In this latter reference a geometric setting was developed replacing the module variety
$mod(A,d)$ for general triangulated categories, mimicking
for this purpose Yoshino's scheme theoretic
approach~\cite{Yoshino}. Various results were given that ensure
that $\leq_\Delta$ or $\leq_{\Delta+\text{nil}}$ define partial
orders on the isomorphism classes of objects of
$\mathcal T$.

Some authors define $\leq_\Delta$ (resp. $\leq_{\Delta+\text{nil}}$) by the
existence of a distinguished triangle $$N\ra M\oplus Z\stackrel{(u_r\;v_r)}{\lra} Z\ra N[1],$$
and some define it as the existence of a distinguished triangle
$$Z\stackrel{{u_\ell}\choose {v_\ell}}{\lra} M\oplus Z\ra N\ra Z[1].$$
Passing to the opposite category the two definitions relate to each other.
Note that the opposite category of a triangulated category is triangulated
as well. However, we show in this paper that actually the situation is even better.
The two possible definitions lead to the same relation on the isomorphism classes of objects
in two important cases. Our main result is the following. In its statement and in the rest
of the paper `artinian ring' means `left and right artinian'.

\begin{Theorem} \label{maintheoremintroduction}
Let $K$ be a commutative ring and
let $\mathcal T$ be a $K$-linear triangulated category  satisfies one of the following two hypotheses
\begin{enumerate}
\item[(a)] Idempotents split in $\mathcal{T}$ and all endomorphism algebras of objects
are artinian,
\item[(b)] $\mathcal T$ is the category of compact objects in an algebraic
compactly generated triangulated $K$-category.
\end{enumerate}
Then for any objects
$M,N$ of $\mathcal{T}$, the following assertions are equivalent:
\begin{enumerate}
\item There is an object $Z_\ell$ of $\mathcal{T}$ and a distinguished triangle in $\mathcal{T}$
$$Z_\ell\stackrel{\begin{pmatrix} v_\ell\\ u_\ell\end{pmatrix}}{\longrightarrow}
Z_\ell\oplus M\stackrel{
}{\longrightarrow} N\longrightarrow
Z_\ell[1],$$
where $v_\ell$ is a nilpotent endomorphism of $Z_\ell$.
\item There is an object $Z_r$ of $\mathcal{T}$ and a distinguished triangle in $\mathcal{T}$
$$N\stackrel{
}{\longrightarrow}
M\oplus Z_r\stackrel{\begin{pmatrix} u_r & v_r\end{pmatrix}}{\longrightarrow} Z_r
\longrightarrow N[1],$$
where $v_r$ is a nilpotent endomorphism of $Z_r$.
\end{enumerate}
\end{Theorem}

It should be noted that since any (left or right) artinian ring is semiperfect (see
\cite[Examples VIII.4]{St})), under the first situation of the theorem, the category
$\mathcal{T}$ is Krull-Schmidt (see \cite[Theorem A.1]{CYZ}).
Moreover, under this hypothesis the assumption
that $v_\ell$ (resp. $v_r$) is nilpotent is inessential.
Indeed, a Fitting lemma type argument can then by applied and this shows that
we can split off a trivial distinguished triangle as direct factor
such that the remaining direct factor distinguished triangle satisfies the nilpotency
hypothesis (cf Remark~\ref{Fittingstyleargument} below
for more details). We cannot avoid the artinian
hypothesis in the first and the nilpotency hypothesis in both cases. The
proof in the first case follows Zwara's arguments in \cite{Zwara}
in  the classical case, but there are quite a few subtleties arising by the
non-uniqueness in the TR3-axiom of triangulated categories.
Zwara frequently uses pushouts and pullbacks and in particular
universal properties which come along with these concepts.
We replace these constructions by homotopy cartesian squares,
and have to cope with the lack of uniqueness of the related construction.
The proof in the second case
is much more involved and heavily uses the concepts developed in \cite{S}.
The main idea in this approach is to use a
dualisation functor like the $K$-duals for ordinary $K$-algebras $A$.
However the situation is more involved here. The hypothesis that the
triangulated category is the category of compact objects in an
algebraic and compactly generated triangulated category
gives that it is actually equivalent to the category of compact objects
in the derived category of some small dg-category.
Then, the new approach is to see this derived category as the
derived category $\mathcal{D}(A)$ of  some dg algebra without unit $A$,
but with sufficiently
many idempotents in a certain sense. Then, it can be shown that one may
dualise with respect to $A$, using the derived functor of the suitable
contravariant $Hom$ functor to $A$.
Further we use in particular the main result of \cite{SZ} in full generality.
The theory of dg algebras with enough idempotents
parallels in a certain sense the development of dg categories
as given by Keller but the situation is new. The approach is presented
in \cite{S}, and we believe that such a theory
is highly useful and should provide many further applications.

The paper is organised as follows. In Section~\ref{dgcategoriessection} we give a
summary of the contents of reference \cite{S}, in order to provide the vocabulary needed
to understand the proof of the main result in the main body of the paper,
without being obliged to go into the full details of that reference. In Section~\ref{recallsection}
we give the relevant background, facts and definitions of degenerations of objects in module categories, as well as in triangulated categories as it was shown in our earlier papers \cite{JSZ,SZ}.
Section~\ref{KrullSchmidtcategorysection} then proves the main result
Theorem~\ref{maintheoremintroduction} under the hypothesis (a), i.e.
in case all objects in the triangulated category have artinian endomorphism ring.
The final Section~\ref{thenonKScasesection}
then gives the proof of
Theorem~\ref{maintheoremintroduction} under hypothesis (b), i.e.
in the case of a triangulated category
which is the category of compact objects in an algebraic compactly generated
triangulated category.

\section{Review on triangulated categories, dg-categories and dg-algebras with enough idempotents}

\label{dgcategoriessection}

For the proof of Theorem~\ref{maintheoremintroduction} under hypothesis (b), which will cover Section~\ref{thenonKScasesection},
we shall need some concepts and statements from the theory of dg-algebras, dg-categories
and triangulated categories in general which are not standard. In particular in case
of categories which do not satisfy Krull-Schmidt theorem, we proceed by considering
dg algebras without units, but having enough idempotents. The complete theory can be found in \cite{S}.
In order to facilitate the reading we summarize  the results of this latter reference and
introduce this way also the notations used in Section~\ref{thenonKScasesection}.
All throughout the rest of the paper,  let $K$
be a commutative ring with unit and all categories which appear all assumed to be $K$-categories.
The unadorned symbol $\otimes$ will stand for the tensor product over $K$.

\subsection{dg categories and dg functors} \label{subsectiondgcategories}
Recall that a \emph{differential graded (dg) $K$-module} is a $\Z$-graded $K$-module $V$ with a
graded endomorphism $d:V\longrightarrow V$ of degree $1$ and square $0$, called the
\emph{differential} (here and all throughout the paper, when the term `differential' is
used to denote a graded map $d$, it will be assumed, without further remark,  that
$d\circ d=0$ and that $d$ is graded and of degree $+1$). We denote by $Dg-K$ or
$\mathcal{C}_{dg}K$ the category of dg $K$-modules. The morphism
space $\text{HOM}_K(V,W)$ in this category is again a dg $K$-module, where the homogeneous
component of degree $n$, denoted $\text{HOM}_K^n(V,W)$, consists of the
homogeneous morphisms of degree $n$. The differential is given by
$d_{Hom}(\alpha)=d_W\circ\alpha-(-1)^{|\alpha|}\alpha\circ d_V$, for any homogeneous morphism
$\alpha\in\text{HOM}_K(V,W)$, where $|\;\cdot\;|$ denotes the degree.

A \emph{dg category} $\mathcal A$ (see \cite{K1} or \cite{K2}) is a category such that
the morphism spaces are dg $K$-modules and the composition map $\text{Hom}_\mathcal{A}(B,C)\otimes\text{Hom}_\mathcal{A}(A,B)\longrightarrow\text{Hom}_\mathcal{A}(A,C)$
satisfies Leibniz rule $d(g\circ f)=d(g)\circ f+(-1)^{|g|}g\circ d(f)$, for all
homogeneous morphisms $f\in\text{Hom}_\mathcal{A}(A,B)$ and $g\in\text{Hom}_\mathcal{A}(B,C)$,
where, abusing notation, we have denoted by $d$ the differential on any of the
appearing $Hom$ spaces. The category $Dg-K$ (denoted by $\mathcal{C}_{dg}K$ in \cite{K2})
is the prototype of a dg category. With any such category, one canonically associates its
\emph{$0$-cycle category} $Z^0\mathcal{A}$ and its \emph{$0$-homology category}
$\mathcal{H}^0\mathcal{A}$. Both of them have the same objects as $\mathcal{A}$,
and as morphisms one puts $\text{Hom}_{Z^0\mathcal{A}}(A,A')=Z^0(\mathcal{A}(A,A'))$
and $\text{Hom}_{H^0\mathcal{A}}(A,A')=H^0(\mathcal{A}(A,A'))$, for all
$A,A'\in\text{Ob}(\mathcal{A})$, the composition of morphisms in both cases
being induced by the composition in $\mathcal{A}$. A \emph{dg functor}
$F:\mathcal{A}\longrightarrow\mathcal{B}$ between dg categories  is just a functor
which preserves the grading and the differential of Hom spaces. Any dg functor
$F:\mathcal{A}\longrightarrow\mathcal{B}$  induces corresponding functors
$F=Z^0F:Z^0\mathcal{A}\longrightarrow Z^0\mathcal{B}$ and
$F:H^0F:H^0\mathcal{A}\longrightarrow H^0\mathcal{B}$.

Associated to  $\mathcal{A}$, there is also the \emph{opposite dg category}
$\mathcal{A}^{op}$ and, given another dg category $\mathcal{B}$, there is a
definition of \emph{tensor product of dg categories} $\mathcal{A}\otimes\mathcal{B}$.
A \emph{homological natural transformation} of dg functors $\tau :F\longrightarrow G$
is a natural transformation such that $\tau_A\in Z^0(\text{Hom}_\mathcal{B}(F(A),G(A)))$,
for any object $A\in\mathcal{A}$. If we have dg functors
$F:\mathcal{A}\longrightarrow\mathcal{B}$ and $G:\mathcal{B}\longrightarrow\mathcal{A}$,
then we have induced dg functors $\mathcal{A}^{op}\otimes\mathcal{B}\longrightarrow Dg-K$,
given by $\text{Hom}_\mathcal{B}(F(?),?)$ and $\text{Hom}_\mathcal{A}(?,G(?))$. A
\emph{dg adjunction} is just an adjunction $(F,G)$ of dg functors such that
the natural isomorphism
$\text{Hom}_\mathcal{B}(F(?),?)\stackrel{\cong}{\longrightarrow}\text{Hom}_\mathcal{A}(?,G(?))$
is a homological natural transformation.
See \cite{K1} and \cite[Section 1]{S} for the details concerning dg categories and dg functors.

\subsection{dg categories and dg algebras with enough idempotents}
Any small $K$-category can be viewed as  an \emph{algebra
with enough idempotents}. The latter is a $K$-algebra $A$ with a
distinguished family $(e_i)_{i\in I}$
of orthogonal idempotents such that $\bigoplus_{i\in I}e_iA=A=\bigoplus Ae_i$.
When such an algebra comes with a grading (as an algebra) such that  the $e_i$
are homogeneous of zero degree,  and with a differential $d:A\longrightarrow A$
such that $d(e_i)=0$, for all $i\in I$, and $d$ satisfies Leibniz rule, then $A$
or the pair $(A,d)$ is called a \emph{differential graded (dg) algebra with
enough idempotents}. It is also shown in \cite{S} that such an algebra may be
viewed as a small dg category with $I$ as set of objects. To any such algebra
$A$ one canonically  associates a (non-small) dg category $Dg-A$, whose objects
are \emph{right dg $A$-modules}. A right dg $A$ module is just a graded right
$A$-module $M$ together with a differential $d_M:M\longrightarrow M$ such that
$d_M(xa)=d_M(x)a+(-1)^{|x|}xd_A(a)$, for all homogeneous elements $x\in M$ and $a\in A$.
Here and in the rest of the paper, unless otherwise specified, all modules are assumed
to be unitary.
That is, we assume that $M=MA$ in our case. We denote by $Gr-A$ the
category with objects the graded right $A$-modules and morphisms the
graded $A$-homomorphisms of degree zero. This category comes with a
canonical equivalence $?[1]:Gr-A\longrightarrow Gr-A$, and we put by
$?[n]:=(?[1])^n$ for each $n\in\mathbb{Z}$. Then,  for each pair $(M,N)$
of right dg $A$-modules, the corresponding space of morphisms in $Dg-A$ is given by
$\text{HOM}_A(M,N)=\bigoplus_{n\in\mathbb{Z}}\text{HOM}_A^n(M,N)$, where
$\text{HOM}_A^n(M,N):=\text{Hom}_{Gr-A}(M,N[n])$. The differential of
$\text{HOM}_A(M,N)$ is the restriction of the differential of $\text{HOM}_K(M,N)$
(see the first paragraph of Section~\ref{subsectiondgcategories}.)
One similarly defines the opposite dg algebra with enough idempotents $A^{op}$
and the tensor product $A\otimes B$ of dg algebras with enough idempotents.
One then defines
the dg category $A-Dg$ of left dg modules and that of dg $A-B-$bimodules,
which are equivalent to $Dg-A^{op}$ and $Dg-(B\otimes A^{op})$, respectively.
This allows to treat the theories of left dg modules or dg bimodules over dg
algebras with enough idempotents just as right dg modules.

\subsection{Stable and derived category of a dg algebra with enough idempotents}
The $0$-cycle
(resp. $0$-homology) category of $Dg-A$ is denoted by $\mathcal{C}(A)$
(resp. $\mathcal{H}(A)$). The category $\mathcal{C}(A)$ is a bicomplete
abelian category, with exact sequences as in $Gr-A$, and, apart from this
abelian structure, it also has a Quillen exact structure, called the
\emph{semi-split exact structure}, where the conflations (=admissible short
exact sequences) are those exact sequences which split in $Gr-A$ (see \cite{B}
for the terminology and main properties of exact categories). With this latter
structure $\mathcal{C}(A)$ is Frobenius, that is, $\mathcal{C}(A)$ has enough
projectives and injectives  and the injective objects coincide with the projective
ones.  The stable category of $\mathcal{C}(A)$, which is then triangulated
(see \cite{H}), is precisely $\mathcal{H}(A)$. This latter (triangulated)
category is called the \emph{homotopy category} of $A$. The class of
\emph{quasi-isomorphisms} in $\mathcal{H}(A)$ (i.e. those morphisms which
induce isomorphisms on homology) is a multiplicative system compatible with
the triangulation in the terminology of Verdier (see \cite{Verdier}). The
localization of $\mathcal{H}(A)$ with respect to the class of quasi-isomorphism,
denoted $\mathcal{D}(A)$, is the \emph{derived category} of $A$. It then has a
unique structure of triangulated category such that the canonical functor
$q:\mathcal{H}(A)\longrightarrow\mathcal{D}(A)$ is triangulated. In \cite{S}
(see, Theorem 3.1 in that reference) it is proved that the theory of dg modules over
dg algebras with enough idempotents and their homotopy and derived categories is
equivalent to the corresponding theory over small dg categories (see \cite{K1}
and \cite{K2} for the details of this latter theory). As a consequence of Keller's
famous theorem (see \cite[Theorem 4.3]{K1}), one gets that any algebraic compactly
generated triangulated category is equivalent to $\mathcal{D}(A)$, for some dg
algebra with enough idempotents $A$ (see \cite[Corollary 6.10]{S}).  Recall that
a triangulated category $\mathcal{T}$ is \emph{algebraic} when it is equivalent
to the stable category of some Frobenius exact category, and that it is called
\emph{compactly generated} when $\mathcal{T}$ has coproducts and there is a set
of compact objects $\mathcal C$
in $\mathcal T$ such that
$\bigcap_{C\in\mathcal{C},n\in\mathbb{Z}}\text{Ker}(\text{Hom}_\mathcal{T}(C[n],?))=0$.
Recall that an object $C$  is \emph{compact} when the functor
$\text{Hom}_\mathcal{T}(C,?):\mathcal{T}\longrightarrow Ab$ preserves arbitrary coproducts.

\subsection{Derived functors}
A (right) dg $A$-module $P$ (resp. $I$) is \emph{homotopically projective}
(resp. \emph{homotopically injective}) when the functor
$\text{HOM}_A(P,?):Dg-A\longrightarrow Dg-K$ (resp. $\text{HOM}_A(?,I):(Dg-A)^{op}\longrightarrow Dg-K$)
preserves acyclic dg modules, something which is equivalent to saying that the induced functor
$\text{Hom}_{\mathcal{H}(A)}(P,?):\mathcal{H}(A)\longrightarrow\text{Mod}-K$
(resp. $\text{Hom}_{\mathcal{H}(A)}(?,I):\mathcal{H}(A)^{op}\longrightarrow\text{Mod}-K$)
vanishes on acyclic dg $A$-modules.
As in the case of small dg categories, the canonical functor
$q_A:\mathcal{H}(A)\longrightarrow\mathcal{D}(A)$ has a left adjoint functor
$\Pi_A :\mathcal{D}(A)\longrightarrow\mathcal{H}(A)$, called the \emph{homotopically
projective resolution functor}, and a right adjoint
$\Upsilon_A:\mathcal{D}(A)\longrightarrow\mathcal{H}(A)$, called the
\emph{homotopically injective resolution functor}, both of which are fully faithful
and triangulated. They are so named because $\text{Im}(\Pi_A)$ (resp. $\text{Im}(\Upsilon_A)$)
consists of homotopically projective (resp. homotopically injective) dg $A$-modules.
The counit $\pi :\Pi_A\circ q_A\longrightarrow 1_{\mathcal{H}(A)}$
(resp. unit $\iota :1_{\mathcal{H}(A)}\longrightarrow q_A\circ\Upsilon_A$) of the
adjunction $(\Pi_A,q_A)$ (resp. $(q_A,\Upsilon_A)$) has the property that $\pi_M$
(resp. $\iota_M$) is a quasi-isomorphism, for each dg module $M$, and it is even an
isomorphism when $M$ is homotopically projective (resp. homotopically injective).
Given  a dg functor $F:Dg-A\longrightarrow Dg-B$  which preserves contractible dg modules,
one defines its \emph{left derived functor} (resp. \emph{right derived functor})
$\mathbb{L}F:\mathcal{D}(A)\longrightarrow\mathcal{D}(B)$ (resp.
$\mathbb{R}F:\mathcal{D}(A)\longrightarrow\mathcal{D}(B)$),  as the composition
$$\mathcal{D}(A)\stackrel{\Pi_A}{\longrightarrow}\mathcal{H}(A)\stackrel{F}{\longrightarrow}
\mathcal{H}(B)\stackrel{q_B}{\longrightarrow}\mathcal{D}(B)$$ (resp. $\mathcal{D}(A)\stackrel{\Upsilon_A}{\longrightarrow}\mathcal{H}(A)
\stackrel{F}{\longrightarrow}\mathcal{H}(B)\stackrel{q_B}{\longrightarrow}\mathcal{D}(B)$).
When the dg functor is contravariant, meaning that $F: (Dg-A)^{op}\longrightarrow Dg-B$ a
dg functor, which preserves contractibility, then we define its
\emph{right derived functor} $\mathbb{R}F$ as the composition  $$\mathcal{D}(A)^{op}\stackrel{\Pi_A^o}{\longrightarrow}\mathcal{H}(A)^{op}
\stackrel{F}{\longrightarrow}\mathcal{H}(B)\stackrel{q_B}{\longrightarrow}\mathcal{D}(B).$$
All these derived functors are triangulated since they are composition
of triangulated functors. If moreover $G:Dg-A\longrightarrow Dg-B$ is another
dg functor as above  and $\tau :F\longrightarrow G$ is a homological natural
transformation of dg functors, then one obtains corresponding natural transformations
of triangulated functors, still denoted the same,
$\tau :\mathbb{L}F\longrightarrow\mathbb{L}G$ and
$\tau :\mathbb{R}F\longrightarrow\mathbb{R}G$ in the covariant case, and just
$\tau :\mathbb{R}F\longrightarrow\mathbb{R}G$ in the contravariant case. Not only
that, but any dg adjunction $(F,G)$ of dg functors gives rise to  a corresponding
triangulated adjunction $(\mathbb{L}F,\mathbb{R}G)$ in the covariant case, and
$((\mathbb{R}F)^o,\mathbb{R}G)$ in the contravariant case (see \cite[Proposition 7.13]{S}).

This somehow classical picture is extended in \cite{S} to dg bifunctors.
Concretely, if $A$, $B$ and $C$ are dg algebras with enough idempotents and
$F:(Dg-A)\otimes (Dg-C)\longrightarrow Dg-B$ is a dg functor which preserves
contractibility on both variables, then one defines
$$\mathbb{L}F:\mathcal{D}(A)\otimes\mathcal{D}(C)
\stackrel{\Pi_A\otimes\Pi_C}{\longrightarrow}\mathcal{H}(A)\otimes\mathcal{H}(C)
\stackrel{H^0F}{\longrightarrow}\mathcal{H}(B)
\stackrel{q_B}{\longrightarrow}\mathcal{D}(B),$$ and
$$\mathbb{R}F:\mathcal{D}(A)\otimes\mathcal{D}(C)
\stackrel{\Upsilon_A\otimes\Upsilon_C}{\longrightarrow}
\mathcal{H}(A)\otimes\mathcal{H}(C)
\stackrel{H^0F}{\longrightarrow}\mathcal{H}(B)
\stackrel{q_B}{\longrightarrow}\mathcal{D}(B).$$
When $F$ is contravariant on the first variable, i.e.
when $F:(Dg-A)^{op}\otimes (Dg-C)\longrightarrow Dg-B$ is a dg functor,
one also defines $$\mathbb{R}F:\mathcal{D}(A)^{op}\otimes\mathcal{D}(C)
\stackrel{\Pi_A^o\otimes\Upsilon_C}{\longrightarrow}\mathcal{H}(A)\otimes\mathcal{H}(C)
\stackrel{H^0F}{\longrightarrow}\mathcal{H}(B)\stackrel{q_B}{\longrightarrow}\mathcal{D}(B).$$
The point is that, under suitable conditions (see \cite[Proposition 7.17]{S} for details),
these later bifunctors are triangulated on each variable.

\subsection{Derived $\text{Hom}$ and $\otimes$ functors}
Given dg algebras with enough idempotents $A$, $B$ and $C$ and dg
bimodules ${}_CM_A$, ${}_BX_A$ and ${}_CU_B$, the dg $K$-modules $\text{HOM}_A(M,X)$
and $U\otimes_BX$ have canonical structures of  dg $B-C-$bimodule and
dg $C-A-$bimodule, respectively, but the first one is non-unitary. This
forces to define the `unitarization'
$$\overline{\text{HOM}}_A(M,X):=B\text{HOM}_A(M,X)C,$$ which is then a (now unitary!)
dg $B-C-$bimodule. It is proved in \cite{S} that the assignments
$(M,X)\rightsquigarrow\overline{\text{HOM}}_A(M,X)$ and
$(U,X)\rightsquigarrow U\otimes_BX$ are the definition on objects of dg functors
\begin{eqnarray*}
\overline{\text{HOM}}_A(?,?):(C-Dg-A)^{op}\otimes (B-Dg-A)&\longrightarrow& B-Dg-C\\
?\otimes_B?:(C-Dg-B)\otimes (B-Dg-A)&\longrightarrow& C-Dg-A.
\end{eqnarray*}
One then puts
\begin{eqnarray*}
?\otimes_B^\mathbb{L}X:=\mathbb{L}(?\otimes_BX):\mathcal{D}(B\otimes C^{op})
&\longrightarrow&\mathcal{D}(A\otimes C^{op}),\\
U\otimes_B^\mathbb{L}?:=\mathbb{L}(U\otimes_B?):\mathcal{D}(A\otimes B^{op})
&\longrightarrow&\mathcal{D}(A\otimes C^{op}),\\
\mathbb{R}\text{Hom}_A(M,?)
:=\mathbb{R}(\overline{\text{HOM}}_A(M,?)):
\mathcal{D}(A\otimes B^{op})&\longrightarrow&\mathcal{D}(C\otimes B^{op})\\
\mathbb{R}\text{Hom}_A(?,X):=\mathbb{R}(\overline{\text{HOM}}_A(?,X)):
\mathcal{D}(A\otimes C^{op})^{op}&\longrightarrow&\mathcal{D}(C\otimes B^{op}).
\end{eqnarray*}
By \cite[Theorems 9.1 and 9.5]{S}, the pairs $$(?\otimes_BX:C-Dg-A\longrightarrow
C-Dg-A,\overline{\text{HOM}}_A(X,?):C-Dg-A\longrightarrow C-Dg-B)$$ and
$$(\overline{\text{HOM}}_{B^{op}}(?,X)^o:B-Dg-C\rightarrow
(C-Dg-A)^{op},\overline{\text{HOM}}_A(?,X):(C-Dg-A)^{op}\rightarrow B-Dg-C)$$
are dg adjunctions and, hence, we get adjunctions of triangulated functors
$$(?\otimes_B^\mathbb{L}X:\mathcal{D}(B\otimes C^{op})\longrightarrow
\mathcal{D}(A\otimes C^{op}),\mathbb{R}\text{Hom}_A(X,?):\mathcal{D}(A\otimes C^{op})
\longrightarrow\mathcal{D}(B\otimes C^{op}))$$ and
$$(\mathbb{R}\text{Hom}_{B^{op}}(?,X)^o:\mathcal{D}(C\otimes B^{op})
\longrightarrow\mathcal{D}(A\otimes C^{op}))^{op},\mathbb{R}\text{Hom}_A(?,X):
\mathcal{D}(A\otimes C^{op}))^{op}\longrightarrow\mathcal{D}(C\otimes B^{op})).$$
By  the previous paragraph, one also defines
$$\mathbb{R}\text{HOM}_A(?,?):=\mathbb{R}(\overline{\text{HOM}}_A(?,?)):
\mathcal{D}(A\otimes C^{op})^{op}\otimes\mathcal{D}(A\otimes B^{op})\longrightarrow
\mathcal{D}(C\otimes B^{op}),$$
which is then a functor which is triangulated
in each variable. Moreover, precise conditions are given in \cite[Corollary 9.7]{S}
to have a natural isomorphisms triangulated functors
$\mathbb{R}\text{HOM}_A(M,?)\cong\mathbb{R}\text{Hom}_A(M,?)$ and $\mathbb{R}\text{HOM}_A(?,X)\cong\mathbb{R}\text{Hom}_A(?,X)$. In particular,
by taking $C=K$ in \cite[Corollary 9.7]{S} and its proof, one gets the following
consequence, which will frequently be used in Section~\ref{thenonKScasesection}:

\begin{Prop} \label{prop.bifunctors-onesidedfunctors}
Let $A$ and $B$ be dg algebras with enough idempotents and let $$\mathbb{R}\text{HOM}_A(?,?):=\mathbb{R}(\overline{\text{HOM}}_A(?,?)):
\mathcal{D}(A^{op})\otimes\mathcal{D}(A\otimes B^{op})\longrightarrow\mathcal{D}(B^{op})$$
be the associated bi-triangulated functor. There are  natural isomorphisms of
triangulated functors, for all dg $B-A-$bimodules $X$ and all right dg $A$-modules $M$:
\begin{enumerate}
\item  $\mathbb{R}\text{HOM}_A(?,X)\cong\mathbb{R}\text{Hom}_A(?,X):
    \mathcal{D}(A)^{op}\longrightarrow\mathcal{D}(B^{op})$.
\item $\mathbb{R}\text{HOM}_A(M,?)\cong\mathbb{R}\text{Hom}_A(\Pi_A(M),?):
\mathcal{D}(A\otimes B^{op})\longrightarrow\mathcal{D}(B^{op})$.
\end{enumerate}
\end{Prop}

On the other hand, when $X={}_AA_A$ is the regular dg bimodule associated
to the dg algebra with enough idempotents $A$, one has that the adjunction $(\mathbb{R}\text{Hom}_{A^{op}}(?,A)^o:\mathcal{D}(A^{op})
\longrightarrow\mathcal{D}(A)^{op},\mathbb{R}\text{Hom}_A(?,A):
\mathcal{D}(A)^{op}\longrightarrow\mathcal{D}(A^{op}))$ gives by
restriction quasi-inverse dualities $\text{per}(A^{op})\stackrel{\cong^o}{\longrightarrow}\text{per}(A)$,
where $\text{per}(A)=\mathcal{D}^c(A)$ (resp. $\text{per}(A^{op})=\mathcal{D}^c(A^{op})$)
is the \emph{right (resp. left) perfect derived category} of $A$, i.e. the full
subcategory of $\mathcal{D}(A)$ (resp. $\mathcal{D}(A^{op})$) consisting of
the compact objects. It is this duality what will allow us to pass from the
left version of degeneration to the right version, and vice versa,  in the
proof of Theorem~\ref{maintheoremintroduction} under hypothesis (b).

\section{Review of degeneration in triangulated categories}

\label{recallsection}

We start to recall from \cite{JSZ,SZ}
a few facts on the concept of degeneration of objects in triangulated categories.

\subsection{The module case}

We first recall a classical result due to Zwara and Riedtmann \cite{Riedtmann1,Zwara1}.
Let $A$ be a $k$-algebra
over an algebraically closed field $k$ and two finite dimensional $A$-modules
$M$ and $N$. Then $N$ belongs to the
closure of the orbit of $M$ if and only if there is a finite dimensional $A$-module $Z_r$
and a short exact sequence
$$0\lra N\stackrel{
}{\longrightarrow}
M\oplus Z_r\stackrel{\begin{pmatrix} u_r & v_r\end{pmatrix}}{\longrightarrow} Z_r
\longrightarrow 0,$$
where $v_r$ is a nilpotent endomorphism of $Z_r$. We say in this case that $M$ degenerates to $N$.
Zwara shows in \cite[Theorem 5]{Zwara} that $M$ degenerates to $N$ if and only if there is
a finite dimensional $A$-module $Z_\ell$ of $\mathcal{T}$ and an exact sequence
triangle
$$0\lra Z_\ell\stackrel{\begin{pmatrix} v_\ell\\ u_\ell\end{pmatrix}}{\longrightarrow}
Z_\ell\oplus M\stackrel{
}{\longrightarrow} N\longrightarrow
0,$$
where $v_\ell$ is a nilpotent endomorphism of $Z_\ell$.

\subsection{Generalising degeneration to triangulated categories}

In our previous work \cite{JSZ,SZ}, and independently by work of Yoshino
\cite{Yoshino} in the case
of stable categories of maximal
Cohen-Macaulay modules of local Gorenstein algebras, the
concept of degeneration for modules was generalised to triangulated categories.
Yoshino discovered
in particular the importance of the hypothesis that the induced endomorphism
of $Z_\ell$ (resp. $Z_r$)
is nilpotent.

\begin{Def}
For two objects $M$ and $N$ of a triangulated category $\mathcal T$
we say that $M$ degenerates to $N$ in the triangle sense and write
$M\leq_{\Delta+nil} N$ if and only if
there is an object $Z_r$ of $\mathcal{T}$ and a distinguished triangle in $\mathcal{T}$
$$Z_\ell\stackrel{\begin{pmatrix} u_\ell\\ v_\ell \end{pmatrix}
}{\longrightarrow}
M\oplus Z_\ell
{\longrightarrow} N
\longrightarrow Z_\ell[1],$$
where $v_\ell$ is a nilpotent endomorphism of $Z_\ell$.
\end{Def}

The main purpose of \cite{SZ} was to define a geometric notion of
degeneration along the lines of
\cite{Yoshino}, and to prove that this notion is equivalent with the
notion of degeneration in the
triangle sense. More precisely we gave the following definition.

\begin{Def}\label{degendatadef}
Let $K$ be a commutative ring and let ${\mathcal C}_K^\circ$ be a
$K$-linear triangulated category with split idempotents.

A degeneration data for ${\mathcal C}_K^\circ$ is given by
\begin{itemize}
\item a triangulated category ${\mathcal C}_K$ with split idempotents and a fully faithful
embedding ${\mathcal C}_K^\circ\longrightarrow{\mathcal C}_K$,
\item
a triangulated category   ${\mathcal C}_V$ with split idempotents
and a full triangulated subcategory ${\mathcal C}_V^\circ$,
\item  triangulated functors $\uar_K^V:{\mathcal C}_K\lra {\mathcal C}_V$
and  $\Phi:{\mathcal C}_V^\circ\ra {\mathcal C}_K$, so that
$({\mathcal C}_K^\circ)\uar_K^V\subseteq {\mathcal C}_V^\circ$, when
we view ${\mathcal C}_K^\circ$ as a full subcategory of ${\mathcal
C}_K$,
\item  a natural transformation
$\text{id}_{{\mathcal C}_V}\stackrel{t}{\lra} \text{id}_{{\mathcal C}_V}$
of triangulated functors
\end{itemize}
These triangulated categories and functors should satisfy the following axioms:
\begin{enumerate}
\item \label{4}
For each object $M$ of ${\mathcal C}_K^\circ$ the morphism
$\Phi(M\uar_K^V)\stackrel{\Phi(t_{M\uar_K^V})}{\lra}\Phi(M\uar_K^V)$
is a split monomorphism in ${\mathcal C}_K$.
\item \label{6}
For all objects $M$ of ${\mathcal C}_K^\circ$ we get
$\Phi(\cone(t_{M\uar_K^V}))\simeq M$.
\end{enumerate}
\end{Def}

Degeneration is then given by the following concept.

\begin{Def}\label{degendef}
Given two objects $M$ and $N$ of ${\mathcal C}_K^\circ$ we say that
$M$ degenerates to $N$ in the categorical sense if there is a
degeneration data for ${\mathcal C}_K^\circ$ and an object $Q$ of
${\mathcal C}_V^\circ$ such that
$$p(Q)\simeq p(M\uar_k^V)\mbox{ in ${\mathcal C}_V^\circ[t^{-1}]$
and }\Phi(\cone(t_Q))\simeq N,$$ where $p:{\mathcal
C}_V^\circ\longrightarrow{\mathcal C}_V^\circ [t^{-1}]$ is the
canonical functor. In this case we write $M\leq_{cdeg}N$.
\end{Def}

\begin{Example}
Yoshino observed that in a triangulated
category $\mathcal T$ for all objects $X$ we get $0\leq_{\Delta+nil} X\oplus X[1]$.
Indeed, $X \ra 0\ra X[1]\stackrel{id}{\ra} X[1]$ and $0\ra X\stackrel{id}{\ra} X\ra 0$
are distinguished triangles for each object $X$. Hence their direct sum
$X\stackrel{0}{\ra} X\oplus 0\ra X\oplus X[1]\ra X[1]$ is a distinguished triangle as well.
Taking $Z=X$ we get the result (in the left version of $\leq_{\Delta+nil}$).

Now, what about the degeneration data interpretation, which is equivalent
to the triangle version in important cases? Then there is an object $Q$ in
some triangulated category ${\mathcal C}_V^\circ$ and an element $t$ in its centre such that
$\Phi(\cone(t_Q))\simeq X\oplus X[1]$ and $p(Q)\simeq 0$ in ${\mathcal C}_V^\circ[t^{-1}]$,
where $p$ is the localisation functor. The latter isomorphism is equivalent to the fact
that $t_Q$ is nilpotent on $Q$. Hence we cannot assume, and actually do not assume, that $Q$
is $t$-flat, as Yoshino does in the case of modules \cite{Yoshinomodules}.
\end{Example}

\subsection{When triangle degeneration is the same as categorical degeneration}

\label{mainresultszsection}

The main result of \cite{SZ} is the following.

\begin{Theorem}\label{mainSZ} Let $K$ be a commutative ring and
let ${\mathcal C}_K^\circ$ be a triangulated $K$-category with split
idempotents. If $M$ and $N$ are objects of $\mathcal{C}_K^\circ$, then
$M\leq_{cdeg}N\Rightarrow M\leq_{\Delta+nil} N$.
When $\mathcal{C}_K^\circ$ is
equivalent to the category of compact objects of a  compactly
generated algebraic triangulated $K$-category, the converse is also
true.
\end{Theorem}

In order to prove that $\leq_{\Delta+nil}$ implies $\leq_{cdeg}$
for the category of compact objects of a  compactly
generated algebraic triangulated $K$-category, we need to construct a degeneration data.
By a result of Keller \cite{K1}, we know
that $\mathcal T$ is equivalent to the category $\mathcal{D}^c(\mathcal{A})$ of
compact objects of the derived category $\mathcal{D}(\mathcal{A})$ of some  small dg-category $\mathcal A$.
We construct then the degeneration data for $\mathcal{C}_K^o=\mathcal{D}^c(\mathcal{A})$ very explicitly,
constructing a dg category $\mathcal{A}[[T]]$
from $\mathcal A$ and our proof of Theorem~\ref{maintheoremintroduction} under hypothesis (b),
which is given in Section~\ref{thenonKScasesection}, uses this construction.
More precisely, recall from the proof of \cite[Proposition 9]{SZ} that if $\mathcal{A}$
is a small dg category, then, considering a variable $T$,  one can form a
new dg category $\mathcal{A}[[T]]$ with the same set of objects as $\mathcal{A}$
and where one defines
$$\text{Hom}_{\mathcal{A}[[T]]}^n(A,A')=
\{\sum_{k\in\mathbb{N}}\alpha_kT^k\text{: }
\alpha_k\in\text{Hom}_\mathcal{A}^n(A,A')\text{, for all }k\in\mathbb{N}\}.$$
Moreover,  one gets a canonical functor
$$?\hat{\otimes}V:\mathcal{C}(\mathcal{A})\longrightarrow\mathcal{C}(\mathcal{A}[[T]])$$
which takes a right dg $\mathcal{A}$-module $M$ to the right dg $\mathcal{A}[[T]]$-module
$M[[T]]:\mathcal{A}^{op}\longrightarrow\mathcal{C}_{dg}K$ acting on objects as
$M[[T]]^n(A)=M^n(A)[[T]]$, for all $n\in\mathbb{Z}$ and all $A\in\mathcal{A}$.
The degeneration data for $\mathcal{C}_K^o$ is then given by taking
$\mathcal{C}_K=\mathcal{D}(\mathcal{A})$, with the corresponding inclusion functor as
$\mathcal{C}_K^o\longrightarrow\mathcal{C}_K$,
$\mathcal{C}_V=\mathcal{D}(\mathcal{A}[[T]])$,
$\mathcal{C}_V^o=\mathcal{D}^c(\mathcal{A}[[T]])$,
the triangulated version of $?\hat{\otimes}V$ as functor $\uar_K^V:\mathcal{C}_K=\mathcal{D}(\mathcal{A})
\longrightarrow\mathcal{D}(\mathcal{A}[[T]])=\mathcal{C}_V$,
the restriction of scalars functor as $\phi :\mathcal{C}_V^o=\mathcal{D}^c(\mathcal{A}[[T]])\longrightarrow\mathcal{D}(\mathcal{A})$
and the natural transformation $t:id_{\mathcal{C}_V}\longrightarrow id_{\mathcal{C}_V}$
is defined by the maps $t_Q:Q\longrightarrow Q$, for each dg
$\mathcal{A}[[T]]$-module $Q$, given by multiplication by $T$.

We take the opportunity to mention that in \cite{SZ} we forgot to mention the grading,
although it was implicit in all of the proofs because, when dealing with dg categories and dg modules,
normally one only uses homogeneous elements. However, a potential reader of \cite{SZ} might think
that there is an error in the definition of $\mathcal{A}[[T]]$ and of $M[[T]]$ for,
as it is written, they are not a graded category or a graded module. The proof of \cite{SZ}
does need not to be modified.

\section{The case of  triangulated categories whose objects have artinian endomorphism algebras}

\label{KrullSchmidtcategorysection}

In this section we shall prove the first part of the theorem.
We will mimic Zwara's proof to give the analogous statement for triangle degeneration.

\subsection{Generalities on homotopy cartesian squares in triangulated categories}

\label{homotopysquaressection}
\label{leftrightequivalent}

As in \cite{Zwara} we first need some preparation. Throughout this section
let $\mathcal T$ be a
triangulated $K$-category for a commutative base ring $K$.
The crucial concept is that of a \emph{homotopy cartesian square}.
Recall from \cite[1.4]{Neeman} that a commutative diagram
$$\xymatrix{
A\ar[r]^a\ar[d]_b&C\ar[d]^c\\B\ar[r]_d&D}$$
is homotopy cartesian if there is a map $e:D\ra A[1]$ such that
$$\xymatrix{A\ar[r]^-{b\choose a}&B\oplus C\ar[r]^-{(-c,d)}&D\ar[r]^e&A[1]}$$
is a distinguished triangle. In the rest of this section, denote by $C_f$ the cone of any morphism $f$ in $\mathcal{T}$. By \cite[Lemma 1.4.3]{Neeman}, if
a commutative square
$$\xymatrix{
A\ar[r]^a\ar[d]_b&C\ar[d]^c\\B\ar[r]_d&D}$$
is homotopy cartesian, then there is an isomorphism $\sigma:C_b\stackrel{\cong}{\longrightarrow}C_d$
such that
$$\xymatrix{
A\ar[r]^a\ar[d]_b&C\ar[d]^c\ar[r]&C_a\ar[d]^\sigma\ar[r]&A[1]\ar[d]^{b[1]}\\
B\ar[r]_d&D\ar[r]&C_d\ar[r]&B[1]}$$
is commutative with rows being distinguished triangles.
\cite[Lemma 1.4.4]{Neeman} gives a partial converse. If
$$\xymatrix{
A\ar[r]^a\ar[d]_b&C\ar[d]^c\ar[r]&C_a\ar[d]^\sigma\ar[r]&A[1]\ar[d]^{b[1]}\\
B\ar[r]_d&D\ar[r]&C_d\ar[r]&B[1]}$$
is commutative, where the rows are distinguished triangles and $\sigma$ is an isomorphism, then
there is a possible different $b':A\ra B$ such that
still
$$\xymatrix{
A\ar[r]^a\ar[d]_{b'}&C\ar[d]^c\ar[r]&C_a\ar[d]^\sigma\ar[r]&A[1]\ar[d]^{b'[1]}\\
B\ar[r]_d&D\ar[r]&C_d\ar[r]&B[1]}$$
is commutative and furthermore
$$\xymatrix{
A\ar[r]^a\ar[d]_{b'}&C\ar[d]^c\\B\ar[r]_d&D}$$
is homotopy cartesian. We may alternatively modify $c$ instead of $b$.

This problem has the annoying consequence that if we have two homotopy cartesian squares
$$\xymatrix{A\ar[r]^a\ar[d]_b&\ar[d]^cC&\mbox{ and }&C\ar[r]^e\ar[d]_c&\ar[d]^fE\\B\ar[r]^d&D
&&D\ar[r]^g&F}
$$
then there is $b'$ and $f'$ fitting in certain morphisms of triangles, such that
$$\xymatrix{A\ar[r]^{ea}\ar[d]_{b'}&\ar[d]^{f'}E\\B\ar[r]^{gd}&F}$$
is a homotopy cartesian square. We would like to be able to assume that $b'=b$ and $f'=f$.
However, we do not know if this is true. Nevertheless, we prove a weaker statement which
satisfy our needs.



\begin{Lemma}\label{pullbacksofpullbacks}
Let $$\xymatrix{X_1\ar[r]^{u_1}\ar[d]^{}&X_2\ar[d]^{f_2}&\text{and}&
X_2\ar[r]^{\nu_1}\ar[d]^{f_2}&X_3\ar[d]^{f_3}\\
0\ar[r]^{}&Y_2&&Y_2\ar[r]^{\nu_2}&Y_3}$$
be homotopy cartesian squares.
Then $$\xymatrix{X_1\ar[r]^{\nu_1u_1}\ar[d]^{}&X_3\ar[d]^{f_3}\\
0\ar[r]^{}&Y_3}$$
is a homotopy cartesian square.
\end{Lemma}

\begin{proof}
We first apply the octahedron axiom to the composition $\nu_1u_1$.
We hence obtain a morphism $v:X_3\ra C_{\nu_1u_1}$, a morphism $\omega:Y_2\ra C_{\nu_1u_1}$,  a morphism $C_{\nu_1u_1}\ra C_{\nu_1}$, and $C_{\nu_1}\ra Y_2[1]$
as indicated below,
$$
\xymatrix{
X_1\ar[r]^{u_1}\ar[drr]_{\nu_1u_1}&X_2\ar[rd]^{\nu_1}\ar[rrr]^{f_2}&&&
Y_2\ar[r]^+\ar@{-->}^{\omega}[dd]&X_1[1]\\
&&X_3\ar[ddrr]\ar^v[drr]&&\\
&&&&C_{\nu_1u_1}\ar[drr]^+\ar@{-->}[d]\\
&&&&C_{\nu_1}\ar[dr]^+\ar@{-->}[d]^+&&X_1[1]\\
&&&&Y_2[1]&X_2[1]
}$$
such that all straight sequences represent distinguished triangles, and such that the diagram is commutative. In particular, $v\nu_1=\omega f_2$. Neeman's interpretation of the octahedral axiom
\cite[Proposition 1.4.6]{Neeman} implies that we may choose $v$ and $\omega$ such that
$$
\xymatrix{X_2\ar[d]_{f_2}\ar[r]^{\nu_1}&X_3\ar[d]^{v}\\ Y_2\ar[r]_\omega&C_{\nu_1u_1} }
$$
is a homotopy cartesian square. Since
$$\xymatrix{X_2\ar[r]^{\nu_1}\ar[d]_{f_2}&X_3\ar[d]^{f_3}\\
Y_2\ar[r]^{\nu_2}&Y_3}$$
is homotopy cartesian by hypothesis, there is an isomorphism $\varphi:C_{\nu_1u_1}\stackrel{\cong}{\longrightarrow} Y_3$
such that
$$\xymatrix{X_2\ar[r]^{\nu_1}\ar[d]_{f_2}&X_3\ar[d]_{v}\ar@/^/[ddr]^{f_3}\\
Y_2\ar[r]^{\omega}\ar@/_/[drr]_{\nu_2}&C_{\nu_1u_1}\ar[dr]^{\varphi}\\
&&Y_3}$$
is commutative. This then shows that
$$
\xymatrix{X_1\ar[r]^{\nu_1u_1}&X_3\ar[r]^{f_3}&Y_3\ar[r]&X_1[1]}
$$
is a distinguished triangle, and hence
$$\xymatrix{X_1\ar[r]^{\nu_1u_1}\ar[d]^{}&X_3\ar[d]^{f_3}\\
0\ar[r]^{}&Y_3}$$
is a homotopy cartesian square as claimed.
\end{proof}

\medskip

A partial converse is true in general however.

\begin{Lemma}\label{inverseofpullbacks}
Let $$\xymatrix{X_1\ar[r]^{u_1}\ar[d]^{f_1}&X_2\ar[d]^{f_2}&&\text{and}&&
X_2\ar[r]^{\nu_1}\ar[d]^{f_2}&X_3\ar[d]^{f_3}\\
Y_1\ar[r]^{u_2}&Y_2&&&&Y_2\ar[r]^{\nu_2}&Y_3}$$
be commutative diagrams. If the first square and
$$\xymatrix{X_1\ar[r]^{\nu_1u_1}\ar[d]^{f_1}&X_3\ar[d]^{f_3}\\
Y_1\ar[r]^{\nu_2u_2}&Y_3}$$
are homotopy cartesian squares, then there is a morphism $\hat\nu_2$
with $\nu_2u_2=\hat\nu_2u_2$ such that
$$\xymatrix{X_2\ar[r]^{\nu_1}\ar[d]^{f_2}&X_3\ar[d]^{f_3}\\
Y_2\ar[r]^{\hat\nu_2}&Y_3}$$ is a homotopy cartesian square.
\end{Lemma}

\begin{proof}
We supposed that
$$\xymatrix{X_1\ar[r]^{u_1}\ar[d]^{f_1}&X_2\ar[d]^{f_2}&\text{and}&
X_1\ar[r]^{\nu_1u_1}\ar[d]^{f_1}&X_3\ar[d]^{f_3}\\
Y_1\ar[r]^{u_2}&Y_2&&Y_1\ar[r]^{\nu_2u_2}&Y_3}$$
are homotopy cartesian squares. Hence we get distinguished triangles

$$\xymatrix{&&X_2\ar[rd]^{-f_2}\\
Y_2[-1]\ar[r]^{\sigma_1[-1]}&X_1\ar[ur]^{u_1}\ar[rd]^{f_1}&\oplus&
Y_2\ar[r]^{\sigma_1}&X_1[1]\\
&&Y_1\ar[ur]_{u_2}
}$$
and

$$\xymatrix{&&X_3\ar[rd]^{-f_3}\\
Y_3[-1]\ar[r]^{\sigma_2[-1]}&X_1\ar[ur]^{\nu_1u_1}\ar[rd]^{f_1}&
\oplus&Y_3\ar[r]^{\sigma_2}&X_1[1]\\
&&Y_1\ar[ur]_{\nu_2u_2}
}$$

We further get a morphism of distinguished triangles by the
fact that the second comes from a square
which factors through the first one.

$$\xymatrix{&&X_2\ar[rd]^{-f_2}\ar@/^1pc/[ddddd]^{\nu_1}\\
Y_2[-1]\ar[r]^{\sigma_1[-1]}\ar[ddddd]_{\nu_2[-1]}&
X_1\ar[ur]_{u_1}\ar[rd]^{f_1}\ar@{=}[ddddd]&\oplus&
Y_2\ar[r]^{\sigma_1}\ar[ddddd]^{\nu_2}&X_1[1]\ar@{=}[ddddd]\\
&&Y_1\ar[ur]_{u_2}\ar@/_1pc/@{=}[ddddd]\\
\\
\\
&&X_3\ar[rd]^{-f_3}\\
Y_3[-1]\ar[r]^{\sigma_2[-1]}&X_1\ar[ur]\ar[rd]_{f_1}&
\oplus&Y_3\ar[r]^{\sigma_2}&X_1[1]\\
&&Y_1\ar[ur]_{\nu_2u_2}
}$$

By Neeman's axiom TR4' we maybe need to modify $\nu_2$ to another map $\hat\nu_2$, forming still a map of
distinguished triangles, so that the cone of this commutative diagram is a distinguished triangle.

$$\xymatrix{&&X_2\ar[rd]^{-f_2}\ar@/^1pc/[ddddd]^{\nu_1}\\
Y_2[-1]\ar[r]^{\sigma_1[-1]}\ar[ddddd]_{\hat\nu_2[-1]}&
X_1\ar[ur]_{u_1}\ar[rd]^{f_1}\ar@{=}[ddddd]&\oplus&
Y_2\ar[r]^{\sigma_1}\ar[ddddd]^{\hat\nu_2}&X_1[1]\ar@{=}[ddddd]\\
&&Y_1\ar[ur]_{u_2}\ar@/_1pc/@{=}[ddddd]\\
\\
\\
&&X_3\ar[rd]_{-f_3}\\
Y_3[-1]\ar[r]^{\sigma_2[-1]}&X_1\ar[ur]\ar[rd]_{f_1}&\oplus&
Y_3\ar[r]^{\sigma_2}&X_1[1]\\
&&Y_1\ar[ur]_{\nu_2u_2}
}$$
This implies in particular that $\nu_2u_2=\hat\nu_2u_2$.
The cone of this has a direct factor isomorphic to

$$\xymatrix{&&X_3\ar[rd]^{-f_3}\\
Y_2[-1]\ar[r]^{\sigma_3[-1]}&X_2\ar[ur]^{\nu_1}\ar[rd]^{f_2}&
\oplus&Y_3\ar[r]^{\sigma_3}&X_2[1]\\
&&Y_2\ar[ur]_{\hat\nu_2}
}$$
which is therefore a distinguished triangle. This proves the statement.
\end{proof}

\subsection{Triangle degeneration is symmetric with respect to left or right existence of the object $Z$}

\label{mainresultsection}

Recall the following result of Xiao-Wu Chen, Yu Ye and Phu Zhang.

\begin{Prop}\cite[Theorem A.1]{CYZ}\label{ChenKS}
An additive category $\mathcal C$ is Krull-Schmidt if and only if any
idempotent splits, and the endomorphism ring of any object of
$\mathcal C$ is semiperfect.
\end{Prop}

Our main result of this section now is the following.

\begin{Theorem}\label{leftZisrightZ}
Let $K$ be a commutative ring and let ${\mathcal T}$ be a $K$-linear
triangulated category with split
idempotents and such that the endomorphism ring of each object is artinian,
and let $M$ and $N$ be two objects. Then
there is an object $Z_r$ and a distinguished triangle
$$N\ra M\oplus Z_r\stackrel{(u_r, v_r)}\lra Z_r\lra N[1]$$
(with nilpotent $v_r$) if and only if
there is an object $Z_\ell$ and a distinguished triangle
$$Z_\ell\stackrel{{u_\ell}\choose {v_\ell}} \lra M\oplus Z_\ell\lra N\lra Z_\ell[1]$$
(with nilpotent $v_\ell$).
\end{Theorem}

\begin{Rem}\label{Fittingstyleargument}
The hypothesis that each object in $\mathcal T$ has artinian endomorphism ring
implies that $\mathcal T$ is Krull-Schmidt, and moreover that
we get Fitting's lemma for $\mathcal T$. In particular,
and endomorphism $\nu:Z\ra Z$ can be decomposed into $$\nu=\left(\begin{array}{cc}\nu'&0\\0&\nu''\end{array}\right):Z=Z'\oplus Z''\ra Z'\oplus Z''=Z$$
and such that $\nu'$ is an automorphism and $\nu''$ is nilpotent. Splitting off the trivial
triangle $Z_\ell'\stackrel{\nu'}{\lra}Z'_\ell\ra 0\ra Z'_\ell[1]$, respectively
$0\ra Z_r'\stackrel{\nu'}{\lra}Z'_r\ra 0[1]$ we may hence assume that $v_\ell$ and $v_r$ are
nilpotent.
\end{Rem}

\begin{proof} (of Theorem~\ref{leftZisrightZ}) Set $N_1:=N$.
Let $Z_\ell$ be an object and let
$$Z_\ell\stackrel{{u_\ell}\choose {v_\ell}} \lra M\oplus Z_\ell\stackrel{(\nu_1,h_1)}\lra N_1\stackrel{\rho_0}\lra Z_\ell[1]$$
be a distinguished triangle with nilpotent endomorphism $v_\ell$.
Then we form the homotopy pushout (cf \cite[Lemma 1.4.3 and 1.4.4]{Neeman})
$$\xymatrix{
Z_\ell\ar[r]^-{{{u_\ell}\choose {v_\ell}}}\ar[d]^{h_1} &M\oplus Z_\ell\ar[r]^-{(\nu_1,h_1)}\ar[d]^{(\nu_2,h_2)}&
N_1\ar[r]^{\rho_0}\ar[d]^{\text{id}}& Z_\ell[1]\ar[d]^{h_1[1]}\\
N_1\ar[r]^{s_1}&N_2\ar[r]^{t_1}&N_1\ar[r]^{\rho_1}&N_1[1]}$$

We get $h_1=t_1h_2$ from the commutativity of the middle square and $h_1[1]\rho_0=\rho_1$
from the right most square. We form the homotopy pushout
$$\xymatrix{Z_\ell\ar[r]^{{{u_\ell}\choose {v_\ell}}}\ar[d]^{h_2} &M\oplus Z_\ell\ar[d]^{(\nu_3,h_3)}\\
N_2\ar[r]^{s_2}&N_3}$$
and obtain a commutative diagram
with homotopy pushouts on the front face and on the back face.
Denote for short $\hat v:={v_\ell\choose u_\ell}$.
$$\xymatrix{
&Z_\ell\ar[rr]^{\hat v}\ar[dd]^{h_2}\ar[dddl]_{h_1}&&
Z_\ell\oplus M\ar[dd]^{(h_3,\nu_3)}\ar[dddl]_{(h_2,\nu_2)}\\
\\
&N_2\ar[rr]^{s_2}\ar[ld]^{t_1}&&N_3\\
N_1\ar[rr]^{s_1}&&N_2
}$$

Since the back face is a homotopy pushout, and since the diagram
$$\xymatrix{
&Z_\ell\ar[rr]^{\hat v}\ar[dd]^{h_2}&&
Z_\ell\oplus M\ar[dd]^{(h_3,\nu_3)}\ar[dddl]_{(h_2,\nu_2)}\\
\\
&N_2\ar[rr]^{s_2}\ar[ld]^{t_1}&&N_3\\
N_1\ar[rr]^{s_1}&&N_2
}$$
is commutative, there is a non-unique map $t_2:N_3\ra N_2$ making the diagram
$$\xymatrix{
&Z_\ell\ar[rr]^{\hat v}\ar[dd]^{h_2}\ar[dddl]_{h_1}&&
Z_\ell\oplus M\ar[dd]^{(h_3,\nu_3)}\ar[dddl]_{(h_2,\nu_2)}\\
\\
&N_2\ar[rr]^{s_2}\ar[ld]^{t_1}&&N_3\ar[ld]^{t_2}\\
N_1\ar[rr]^{s_1}&&N_2
}$$
commutative. Lemma~\ref{inverseofpullbacks} then allows to modify $t_2$ such
that the bottom face of the diagram is homotopy cartesian, such that all already
shown commutativity properties still hold, and such that the above diagram with modified $t_2$
is still commutative. Here, in order to simplify the notation,
we denote the modified $t_2$ again by $t_2$.

\medskip

We proceed now by induction on the degree $n$. Suppose we have a commutative diagram with rows being
distinguished triangles and whose square faces are homotopy cartesian squares
$$\xymatrix{
&Z_\ell\ar[rr]^{\hat v}\ar[dd]^{h_{n-1}}\ar[dddl]_{h_{n-2}}&&
Z_\ell\oplus M\ar[dd]^{(h_n,\nu_n)}\ar[dddl]_{(h_{n-1},\nu_{n-1})}\\
\\
&N_{n-1}\ar[rr]^{s_{n-1}}\ar[ld]^{t_{n-2}}&&
N_n\ar[ld]^{t_{n-1}}\\
N_{n-2}\ar[rr]^{s_{n-2}}&&N_{n-1}
}$$

We form the homotopy pushout (on the back face of the diagram defining $N_{n+1}$, $h_{n+1}$, $s_n$ and $\nu_{n+1}$)
$$\xymatrix{
&Z_\ell\ar[rr]^{\hat v}\ar[dd]^{h_{n}}\ar[dddl]_{h_{n-1}}&&
Z_\ell\oplus M\ar[dd]^{(h_{n+1},\nu_{n+1})}\ar[dddl]_{(h_n,\nu_n)}\\
\\
&N_{n}\ar[rr]^{s_{n}}\ar[ld]^{t_{n-1}}&&
N_{n+1}\\
N_{n-1}\ar[rr]^{s_{n-1}}&&N_{n}
}$$
Since the diagram
$$\xymatrix{
&Z_\ell\ar[rr]^{\hat v}\ar[dd]^{h_{n}}
&&
Z_\ell\oplus M\ar[dd]^{(h_{n+1},\nu_{n+1})}\ar[dddl]_{(h_n,\nu_n)}\\
\\
&N_{n}\ar[rr]^{s_{n}}\ar[ld]^{t_{n-1}}&&
N_{n+1}\\
N_{n-1}\ar[rr]^{s_{n-1}}&&N_{n}
}$$
is commutative, there is a morphism $t_n:N_{n+1}\ra N_n$ making the diagram
$$\xymatrix{
&Z_\ell\ar[rr]^{\hat v}\ar[dd]^{h_{n}}\ar[dddl]_{h_{n-1}}&&
Z_\ell\oplus M
\ar[dd]^{(h_{n+1},\nu_{n+1})}\ar[dddl]_{(h_n,\nu_n)}\\
\\
&N_{n}\ar[rr]^{s_{n}}\ar[ld]^{t_{n-1}}&&
N_{n+1}\ar[dl]^{t_n}\\
N_{n-1}\ar[rr]^{s_{n-1}}&&N_{n}
}$$
commutative. By Lemma~\ref{inverseofpullbacks} we may modify $t_n$ without
changing commutativity of what is already shown (we denote the modified
$t_n$ again by $t_n$), such that the bottom face of the diagram is a
homotopy pushout and such the above diagram
is still commutative.

We now continue as in the proof of \cite[Theorem 2.3]{Zwara}.
We define $\omega_1:=\nu_1$ and
$$\omega_{j+1}:=(\nu_{j+1},s_j\omega_j)\in
Hom_{\mathcal T}(M\oplus M^{j},N_{j+1})
=Hom_{\mathcal T}(M^{j+1},N_{j+1})$$
for all $j\geq 1$.

In the diagram
$$\xymatrix{
Z_\ell\ar[r]^{{{u_\ell}\choose {v_\ell}}}\ar[d]^{h_j}
&M\oplus Z_\ell\ar[rr]^{(\nu_1,h_1)}\ar[d]^{(\nu_{j+1},h_{j+1})}&&
N\ar[r]^{\rho}\ar[d]^{\text{id}}& Z_\ell[1]\ar[d]^{h_2[1]}\\
N_j\ar[r]^{s_j}&N_{j+1}\ar[rr]^{v_j}&&N\ar[r]^{\rho_j}&N_j[1]
}$$
the left square is a homotopy cartesian square. Therefore we get a distinguished triangle
$$\xymatrix{
Z_\ell\ar[rr]^-{\left(\begin{array}{c}u_\ell\\v_\ell\\h_j\end{array}\right)}&&
M\oplus Z_\ell\oplus N_j\ar[rrr]^-{(\nu_{j+1},h_{j+1},-s_j)}&&&N_{j+1}\ar[r]^{\sigma_j}&Z_\ell[1]
}$$

Since
$$\xymatrix{M^j\ar[r]^{\text{id}}&M^j\ar[r]&0\ar[r]&M^{j}[1]}$$
is a distinguished triangle, also
the direct sum of these two distinguished triangles
$$\xymatrix{
M^j\oplus Z_\ell\ar[rr]^-{\left(
\begin{array}{cc}1_{M^j}&0\\0&u_\ell\\0&v_\ell\\0&h_j\end{array}\right)}&&
M^j\oplus M\oplus Z_\ell\oplus N_j\ar[rrr]^-{(0,\nu_{j+1},h_{j+1},-s_j)}&&&
N_{j+1}\ar[r]^-{(0,\sigma_j)}&(M^j\oplus Z_\ell)[1]
}$$
is distinguished.
For the series of morphisms $\omega_j:M^j\ra N_j$ satisfying
$\omega_1=\nu_1$ and $\omega_{j+1}=(\nu_{j+1},s_j\omega_j)$
for all $j$, this distinguished triangle
is isomorphic to the triangle
$$\xymatrix{
M^j\oplus Z_\ell\ar[rr]^-{\left(
\begin{array}{cc}1_{M^j}&0\\0&u_\ell\\0&v_\ell\\\omega_j&h_j\end{array}\right)}&&
M^j\oplus M\oplus Z_\ell\oplus N_j\ar[rrr]^-{(\omega_{j+1},\nu_{j+1},h_{j+1},-s_j)}&&&
N_{j+1}\ar[r]^-{(0,\sigma_j)}&(M^j\oplus Z_\ell)[1]
}$$
Indeed, we get morphisms of triangles
$$\xymatrix{M^j\ar[rr]^{\text{id}}&&M^j\ar[rrr]&&&0\ar[r]&M^j[1]\\ \\
M^j\oplus Z_\ell\ar[uu]^{(1,0)}\ar[rr]^-{\left(
\begin{array}{cc}1_{M^j}&0\\0&u_\ell\\0&v_\ell\\\omega_j&h_j\end{array}\right)}&&
M^j\oplus M\oplus Z_\ell\oplus N_j\ar[rrr]^-{(\omega_{j+1},\nu_{j+1},h_{j+1},-s_j)}
\ar[uu]_{(1,0,0,0)}&&&N_{j+1}\ar[r]\ar[uu]&(M^j\oplus Z_\ell)[1]\ar[uu]\\ \\ \\
M^j\ar[rr]^{\text{id}}
\ar[uuu]^{\left(\begin{array}c 1\\0\end{array}\right)}&&
M^j\ar[rrr]\ar[uuu]^{\left(
\begin{array}c 1\\0\\0\\0\\\omega_j\end{array}\right)}&&&0\ar[r]\ar[uuu]&M^j[1]\ar[uuu]
}$$
and therefore the middle triangle has a direct factor
$$\xymatrix{M^j\ar[r]^{\text{id}}&M^j\ar[r]&0\ar[r]&M^{j}[1]}$$
and the remaining direct factor is the original distinguished triangle.

Hence,
$$\xymatrix{
M^j\oplus Z_\ell
\ar[rr]^{\scriptsize\left(\begin{array}{cc}1_{M^j}&0\\0&u_\ell\\0&v_\ell
\end{array}\right)}
\ar[d]_{(\omega_j,h_j)}&&M^{j+1}\oplus Z_\ell\ar[d]^{(\omega_{j+1},h_{j+1})}\\
N_j\ar[rr]^{s_j}&&N_{j+1}
}$$
is a homotopy cartesian square. Moreover,
$$\xymatrix{
Z_\ell\ar[r]^{u_\ell\choose v_\ell}\ar[d]&M\oplus Z_\ell\ar[d]^{(\omega_1,h_1)}\\ 0\ar[r]&N
}$$
is a homotopy cartesian square by hypothesis.
Now, since $$\xymatrix{
Z_\ell\ar[r]^{u_\ell\choose v_\ell}\ar[d]&M\oplus Z_\ell\ar[d]^{(\omega_1,h_1)}&\mbox{ and }
&M^j\oplus Z_\ell
\ar[rr]^{\scriptsize\left(\begin{array}{cc}1_{M^j}&0\\0&u_\ell\\0&v_\ell
\end{array}\right)}
\ar[d]_{(\omega_j,h_j)}&&M^{j+1}\oplus Z_\ell\ar[d]^{(\omega_{j+1},h_{j+1})}\\
0\ar[r]&N&&N_j\ar[rr]^{s_j}&&N_{j+1}
}$$
are homotopy cartesian squares, applying
Lemma~\ref{pullbacksofpullbacks},
we get a homotopy cartesian square
$$\xymatrix{
Z_\ell\ar[r]^{\psi_2\choose v_\ell^2}\ar[d]&M^2\oplus Z_\ell\ar[d]^{(\omega_2,h_2)}\\ 0\ar[r]&N_2
}$$
and by induction on $j$ we get that there
is a morphism $\psi_j:Z_\ell\ra M^j$
and a homotopy cartesian square
$$\xymatrix{
Z_\ell\ar[r]^{\psi_j\choose v_\ell^j}\ar[d]&M^j\oplus Z_\ell\ar[d]^{(\omega_j,h_j)}\\ 0\ar[r]&N_j
}$$
for all $j\in\N$.
Hence for all $j>0$ there
is a morphism $\psi_j:Z_\ell\ra M^j$ such that
$$Z_\ell\stackrel{\psi_j\choose v_\ell^j}{\lra}M^j\oplus
Z_\ell\stackrel{(\omega_j,h_j)}{\lra}N_j\lra Z_\ell[1]$$
is a distinguished triangle.
Now $v_\ell$ is nilpotent of degree $k_0$, say.
For any $k\geq k_0$ consider the commutative diagram with distinguished triangles
in the horizontal rows
$$\xymatrix{Z_\ell\ar[r]^{\psi_k\choose {v_\ell^k}}\ar[d]&M^k\oplus
Z_\ell\ar[rr]^{(\omega_k,h_k)}\ar[d]^{(0,\text{id})}&&N_k\ar[r]&Z_\ell[1]\ar[d]\\
0\ar[r]&Z_\ell\ar[rr]^{\text{id}}&&Z_\ell\ar[r]&0
}$$
which can be completed by a map $z_k$ to a morphism of distinguished
triangles, using TR3 (see e.g. \cite[Definition 1.1.1]{Neeman} or \cite[Definition 3.4.1]{reptheobuch}),
$$\xymatrix{Z_\ell\ar[r]^{\psi_k\choose {v_\ell^k}}\ar[d]&
M^k\oplus Z_\ell\ar[rr]^{(\omega_k,h_k)}\ar[d]^{(0,\text{id})}&&N_k\ar[r]\ar[d]^{z_k}&
Z_\ell[1]\ar[d]\\
0\ar[r]&Z_\ell\ar[rr]^{\text{id}}&&Z_\ell\ar[r]&0
}$$
Therefore $z_kh_k=\text{id}_{Z_\ell}$, and $h_k$ is a split monomorphism.
We recall that we constructed the sequence $h_j$ as iterated homotopy pushouts, and hence, by
definition we have a homotopy cartesian square
$$\xymatrix{
Z_\ell\ar[r]^-{\hat v}\ar[d]_{h_k}&Z_\ell\oplus M\ar[d]^{(h_{k+1},\nu_{k+1})}\\
N_k\ar[r]_-{s_n}&N_{k+1}
}$$
where $h_k$ (and $h_{k+1}$) are split monomorphisms.
This shows first (cf e.g. \cite[Lemma 3.4.9]{reptheobuch}) that
$N_k\simeq Z_\ell\oplus C_{h_k}$, for $C_{h_k}$ being the mapping cone of $h_k$.
Moreover, we get a distinguished triangle
$$
\xymatrix{Z_\ell\ar[r]^-{\hat v\choose h_k}&Z_\ell\oplus M\oplus
N_k\ar[rrr]^-{(h_{k+1},\nu_{k+1},-s_k)}&&&N_{k+1}\ar[r]&Z_\ell[1]}.
$$
Since $h_k$ is split monomorphism, this distinguished triangle is isomorphic to the direct sum
of the trivial distinguished triangle
$$\xymatrix{Z_\ell\ar[r]^{\text{id}}&Z_\ell\ar[r]&0\ar[r]&Z_\ell[1]}$$
and
$$\xymatrix{0\ar[r]&0\oplus C_{h_k}\oplus Z_\ell\oplus M\ar[r]&N_{k+1}\ar[r]&0[1].}$$
Hence (cf e.g. \cite[Lemma 3.4.9]{reptheobuch}),
$$(\dagger)\;\;N_{k+1}\simeq C_{h_k}\oplus Z_\ell\oplus M\simeq N_k\oplus M.$$
Note that we could have argued also that since $z_k$ is left inverse to $h_k$,
$(0,z_k)$ is left inverse to $\hat v\choose h_k$, and so the above triangle splits.
This then gives the desired isomorphism $(\dagger)$ via Remark~\ref{Fittingstyleargument}.

Recall that, posing $N_0:=0$, for all $j\geq 1$ we have by
construction homotopy cartesian squares
$$\xymatrix{N_j\ar[r]^{s_j}\ar[d]_{t_j}&N_{j+1}\ar[d]^{t_{j+1}}\\
N_{j-1}\ar[r]^{s_j}&N_{j}}$$
for all $j\geq 1$. Using Lemma~\ref{pullbacksofpullbacks} and an obvious induction as above
this implies that we get
an cartesian square
$$\xymatrix{N\ar[r]\ar[d]&N_{k+1}\ar[d]\\ 0\ar[r]&N_k}$$
which gives a distinguished triangle
$$\xymatrix{N\ar[r]&N_{k+1}\ar[r]&N_k\ar[r]&N[1]}$$
where we have chosen $k$ such that $v_\ell^k=0$. By equation
$(\dagger)$ we get $N_{k+1}\simeq N_k\oplus M$, which shows
that there is a distinguished triangle
$$N\lra  N_k\oplus M\lra N_k\lra N[1].$$
Posing $N_k=:Z_r$ this gives the statement, except that we do not
get yet that the induced endomorphism of $Z_r$ is nilpotent.

Since the endomorphism ring of all objects in $\mathcal T$ are artinian and
idempotents split, $\mathcal T$ is Krull-Schmidt (cf Proposition~\ref{ChenKS}), and then we may
split off $f'$ in a nilpotent endomorphism of a direct factor and an
automorphism of a direct factor, using Fitting's lemma. The automorphism
part splits in the distinguished triangle, and we obtain the statement.

The other direction is done applying the statement proved above to the
opposite category ${\mathcal T}^{op}$ of $\mathcal T$.
\end{proof}


\begin{Rem} A triangulated category with split idempotents for which each
object has artinian endomorphism rings is Krull-Schmidt
and a Fitting-like theorem holds (cf Remark~\ref{Fittingstyleargument}).
However, if $\mathcal T$ is a general triangulated category, and
in particular if we do not assume that $\mathcal T$ is Krull-Schmidt,
then we get a weaker statement in Theorem~\ref{leftZisrightZ}.
The hypothesis that $\mathcal T$ has artinian endomorphism
rings is only used at the very end of the proof of the theorem in order to be able to split off a
direct factor in order to get a nilpotent endomorphism on the remaining factor.
If $\mathcal T$ is a general triangulated category
we proved that the existence of a distinguished triangle
$$Z_\ell\stackrel{{u_\ell}\choose {v_\ell}} \lra M\oplus Z_\ell\lra N\lra Z_\ell[1],$$
with nilpotent $v_\ell$ implies the existence of a distinguished triangle
$$N\lra M\oplus Z_r\stackrel{(u_r\;v_r)}\lra Z_r[1]\lra N[1],$$
but we are unable to deduce that $v_r$ is nilpotent.
\end{Rem}

\section{The case of a category of compact objects in an algebraic compactly generated category}

\label{thenonKScasesection}

Recall from Section~\ref{mainresultszsection} the construction of the category
$\mathcal{A}[[T]]$.
When we have a dg algebra with enough idempotents $A$ and view it as a
dg category (see \cite[Section 3]{S}),
we can undertake the same construction, but there is a classical notion of power
series algebra $A[[T]]$ which does not correspond to  the `power series dg category'
mentioned in Section~\ref{mainresultszsection}. The corresponding dg algebra with enough
idempotents is the subalgebra of $A[[T]]$ given as
$$\tilde{A}[[T]]=\bigoplus_{i,j\in I}\bigoplus_{n\in\mathbb{Z}}e_iA^ne_j[[T]].$$
That is,  $\tilde{A}[[T]]$ consists of the power series $\sum_{k\in\mathbb{N}}a_kT^k$
for which there are finite subsets $J\subset I$ and $F\subset\mathbb{Z}$ such that
$a_k\in\bigoplus_{i,j\in J}\bigoplus_{n\in F}e_iA^ne_j$, for all $k\in\mathbb{N}$.
The subalgebra $\tilde{A}[[T]]$ is made into a dg algebra by defining
$\tilde{A}[[T]]^n=\bigoplus_{i,j\in I}e_iA^ne_j[[T]]$, for each $n\in\mathbb{Z}$,
and by defining its differential by the rule
$d(\sum_{k\in\mathbb{N}}a_kT^k)=\sum_{i\in\mathbb{N}}d(a_k)T^k$.
Note that then $(e_i)_{i\in I}$ is also a distinguished family of
orthogonal idempotents of $\tilde{A}[[T]]$. Moreover, we have  a  canonical
inclusion $\iota: A\hookrightarrow\tilde{A}[[T]]$ ($a\rightsquigarrow a=aT^o$)
and a canonical augmentation map
$\rho :\tilde{A}[[T]]\longrightarrow A$ ($\sum_{k\in\mathbb{N}}a_kT^k\rightsquigarrow a_0$)
such that $\rho\circ\iota =1_A$. Both $\iota$ and $\rho$ are
homomorphisms of dg algebras making the codomain into a unitary bimodule over the domain.
We can apply to $\iota$ and $\rho$ the results of  \cite{S}
(see Corollary 9.4 and Section 10 in that reference) concerning homomorphisms
of dg algebras with enough idempotents.

Likewise for right $A$-modules, there is a classical notion of
`module of power series' $M[[T]]$ which does not correspond the the `dg
$\mathcal{A}[[T]]$-module of power series' mentioned in Section~\ref{mainresultszsection}.
The corresponding right dg $\tilde{A}[[T]]$-module is
$$\tilde{M}[[T]]=\bigoplus_{i\in I}\bigoplus_{n\in\mathbb{Z}}M^ne_i[[T]],$$
which is the
$K$-submodule of $M[[T]]$ consisting of those power series $\sum_{k\in\mathbb{N}}m_kT^k$
for which there exist finite subsets $J\subset I$ and $F\subset\mathbb{Z}$,
both depending on the power series, such that
$m_k\in\oplus_{i\in J,n\in F}M^ne_i$, for all $k\in\mathbb{N}$. Defining the grading
by the rule $\tilde{M}[[T]]^n=\bigoplus_{i\in I}M^ne_i[[T]]$ and the differential
$d:\tilde{M}[[T]]\longrightarrow\tilde{M}[[T]]$ by the rule
$d(\sum_{k\in\mathbb{N}}m_kT^k)=\sum_{k\in\mathbb{N}}d(m_k)T^k$, we clearly
endow $\tilde{M}[[T]]$ with a structure of right dg $\tilde{A}[[T]]$-module.
We warn the reader of the possible confusion of this construction when applied
to the regular right dg $A$-module $A_A$ for the resulting right dg $\tilde{A}[[T]]$-module
is not equal to the regular right dg module $\tilde{A}[[T]]_{\tilde{A}[[T]]}$.

\subsection{Dualising degeneration data}

We now give, and actually extend, the functor corresponding to $?\hat{\otimes}V$ to the
context of dg algebras with enough idempotents and dg modules over them.

\begin{Prop} \label{prop.the tilde functor}
The assignment $M\rightsquigarrow\tilde{M}[[T]]$ is the definition on objects
of a dg functor $?\hat{\otimes}V:Dg-A\longrightarrow Dg-\tilde{A}[[T]]$ which
satisfies the following properties:

\begin{enumerate}
\item \label{prop.the tilde functor.1} $?\hat{\otimes}V$ takes contractible dg modules to contractible dg modules.
\item \label{prop.the tilde functor.2} The associated functor on $0$-cycle categories
$$?\hat{\otimes}V:Z^0(Dg-A)=\mathcal{C}(A)\longrightarrow
\mathcal{C}(\tilde{A}[[T]])=Z^0(Dg-\tilde{A}[[T]])$$
is exact with respect to the respective abelian structures.
\item \label{prop.the tilde functor.3} $?\hat{\otimes}V$ preserves acyclic dg modules. In particular,
it induces a triangulated functor $?\hat{\otimes}V:\mathcal{D}(A)\longrightarrow\mathcal{D}(\tilde{A}[[T]])$
which is both the left and right derived functor of `itself'. We will denote this functor by $\uar_K^V$.
\item \label{prop.the tilde functor.4} The multiplication map $M\otimes_A\tilde{A}[[T]]\longrightarrow\tilde{M}[[T]]$
defines a homological natural transformation of dg functors $\mu :\iota^*\longrightarrow ?\hat{\otimes}V$
\item \label{prop.the tilde functor.5} The induced natural transformation $\mu :\mathbb{L}\iota^*\longrightarrow\uar_K^V$ of triangulated functors
$\mathcal{D}(A)\longrightarrow\mathcal{D}(\tilde{A}[[T]])$  is a natural isomorphism when
evaluated at objects of $\text{per}(A)$.
\end{enumerate}
\end{Prop}

\begin{proof} If $f:M\longrightarrow N$ is a homogeneous morphism in $Dg-A$,
we define $\tilde{f}:=(?\hat{\otimes}V)(f)$ by the rule
$\tilde{f}(\sum_{k\in\mathbb{N}}m_kT^k)=\sum_{k\in\mathbb{N}}f(m_k)T^k$.
It is routine to verify that the assignments $M\rightsquigarrow\tilde{M}[[T]]$ and
$f\rightsquigarrow\tilde{f}$ give a graded functor $GR-A\longrightarrow GR-\tilde{A}[[T]]$.
In order to check that they define a dg functor
$?\hat{\otimes}V:Dg-A\longrightarrow Dg-\tilde{A}[[T]]$, we  need to check that
if $M,N$ are right dg $A$-modules and $d:\text{HOM}_{A}(M,N)\longrightarrow\text{HOM}_{A}(M,N)$
and $\delta:\text{HOM}_{\tilde{A}[[T]]}(\tilde{M}[[T]],\tilde{N}[[T]])
\longrightarrow\text{HOM}_{\tilde{A}[[T]]}(\tilde{M}[[T]],\tilde{N}[[T]])$
are the respective differentials on Hom spaces, then one has
$\delta (\tilde{f})=\widetilde{d(f)}$, for any homogeneous element $f\in\text{HOM}_A(M,N)$.
On one hand, we have that
$$\delta (\tilde{f})=d_{\tilde{N}[[T]]}\circ\tilde{f}-(-1)^{|f|}\tilde{f}\circ d_{\tilde{M}[[T]]}. $$
On the other hand, if we let act $\widetilde{d(f)}$ on a homogeneous element
$\sum_{k\in\mathbb{N}}m_kT^k\in\tilde{M}[[T]]$ (whence the degree $deg(m_k)$ is independent of $k$),
then we get:

\begin{eqnarray*}
\widetilde{d(f)}\left(\sum_{k\in\mathbb{N}}m_kT^k\right)&=&\sum_{k\in\mathbb{N}}d(f)(m_k)T^k\\
&=&\sum_{k\in\mathbb{N}}\left[(d_N\circ f-(-1)^{|f|}f\circ d_M)(m_k)\right]T^k\\
&=&\sum_{k\in\mathbb{N}}d_N(f(m_k))T^k-(-1)^{|f|}\sum_{k\in\mathbb{N}}f(d_M(m_k))T^k\\
&=&d_{\tilde{N}[[T]]}\left(\sum_{k\in\mathbb{N}}f(m_k)T^k\right)-(-1)^{|f|}\tilde{f}\left(\sum_{k\in\mathbb{N}}d_M(m_k)T^k\right)\\
&=&\left[d_{\tilde{N}[[T]]}\circ\tilde{f}-(-1)^{|f|}\tilde{f}\circ d_{\tilde{M}[[T]]}\right]\left(\sum_{k\in\mathbb{N}}m_kT^k\right).
\end{eqnarray*}
This shows that  $\delta (\tilde{f})=\widetilde{d(f)}$, as desired.

Finally, it is also routine to see that
$(?\hat{\otimes}V)(\text{cone}(1_M))\cong\text{cone}(1_{(?\hat{\otimes}V)(M)})$, which ends
the proof of assertion (\ref{prop.the tilde functor.1}).

\vspace*{0.3cm}

(\ref{prop.the tilde functor.2}) Let $L\stackrel{f}{\longrightarrow} M\stackrel{g}{\longrightarrow} N$ be an exact
sequence in $\mathcal{C}(A)$, when this category is considered with its natural abelian
structure, and consider the corresponding sequence $\tilde{L}[[T]]\stackrel{\tilde{f}}{\longrightarrow}\tilde{M}[[T]]
\stackrel{\tilde{g}}{\longrightarrow}\tilde{N}[[T]]$.
Since  $\tilde{f}$ and $\tilde{g}$ are both morphisms in $Gr-\tilde{A}[[T]]$ we just
need to check that the induced sequence
$$\tilde{L}[[T]]^ne_i=L^ne_i[[T]]
\stackrel{\tilde{f}}{\longrightarrow}\tilde{M}[[T]]^ne_i=
M^ne_i[[T]]\stackrel{\tilde{g}}{\longrightarrow}\tilde{N}[[T]]^ne_i=N^ne_i[[T]]$$
is exact, for all $i\in I$ and $n\in\mathbb{Z}$.
But, given  $\sum_{k\in\mathbb{N}}m_kT^k\in M^ne_i[[T]]$, we have that
$\tilde{g}(\sum_{k\in\mathbb{N}}m_kT^k)=0$ if and only if $g(m_k)=0$ for all
$k\in\mathbb{N}$. This in turn is equivalent to saying that, for each $k\in\mathbb{N}$,
there exists a $l_k\in L^ne_i$ such that $f(l_k)=m_k$. That is, we have that
$\sum_{k\in\mathbb{N}}m_kT^k\in\text{Ker}(\tilde{g})$ if and only if $\sum_{k\in\mathbb{N}}m_kT^k=\tilde{f}(\sum_{k\in\mathbb{N}}l_kT^k)$,
for some $\sum_{k\in\mathbb{N}}l_kT^k\in L^ne_i[[T]]$.

\vspace*{0.3cm}

(\ref{prop.the tilde functor.3}) Let $M$ be an acyclic right dg $A$-module, let $\sum_{k\in\mathbb{N}}m_kT^k$ be an element of $\text{Ker}(d^n:\tilde{M}[[T]]^n\longrightarrow\tilde{M}[[T]]^{n+1})$ and let $J\subset I$
be a finite subset such that $m_k\in\bigoplus_{i\in J}M^ne_i$, for all $k\in\mathbb{N}$.
By the acyclicity condition of $M$, for each $k\in\mathbb{N}$, we have an
$m'_k\in\bigoplus_{i\in J}M^{n-1}e_i$ such that $d(m'_k)=m_k$. It follows
$\sum_{k\in\mathbb{N}}m'_kT^k$ is an element of $\tilde{M}[[T]]^n$ such that
$d(\sum_{k\in\mathbb{N}}m'_kT^k)=\sum_{k\in\mathbb{N}}m_kT^k$. Therefore
$\tilde{M}[[T]]$ is an acyclic right dg $\tilde{A}[[T]]$-module.
The last comment of the assertion follows from \cite[Remark 7.11]{S}.

\vspace*{0.3cm}

(\ref{prop.the tilde functor.4})  We clearly have that $\mu_M:M\otimes_A\tilde{A}[[T]]\longrightarrow\tilde{M}[[T]]$
is a morphism (of zero degree) in $GR-\tilde{A}[[T]]$. In addition, if we denote by
$d:M\otimes_A\tilde{A}[[T]]\longrightarrow M\otimes_A\tilde{A}[[T]]$ and
$\tilde{d}:\tilde{M}[[T]]\longrightarrow\tilde{M}[[T]]$ the respective differentials,
then, for all homogeneous elements $m\in M$ and $\sum_{k\in\mathbb{N}}a_kT^k$, we have
\begin{eqnarray*}
(\mu_M\circ d)\left[m\otimes (\sum_{k\in\mathbb{N}}a_kT^k)\right]&=&
\mu_M \left(d_M(m)\otimes (\sum_{k\in\mathbb{N}}a_kT^k)+(-1)^{|m|}m\otimes (\sum_{k\in\mathbb{N}}d(a_k)T^k)\right)\\
&=&\sum_{k\in\mathbb{N}}d_M(m)a_kT^k +(-1)^{|m|}md(a_k)T^k\\
&=&\sum_{k\in\mathbb{N}}(d_M(m)a_k+(-1)^{|m|}md(a_k))T^k\\
&=&\sum_{k\in\mathbb{N}}d_M(ma_k)T^k\\
&=&\tilde{d}(\sum_{k\in\mathbb{N}}ma_kT^k)\\
&=&(\tilde{d}\circ\mu_M)\left[m\otimes (\sum_{k\in\mathbb{N}}a_kT^k)\right].
\end{eqnarray*}
Then, once the naturality $\mu$ is proved, we will have that it is actually
a homological natural transformation of dg functors (see \cite[Remark 7.1]{S}). But that naturality is clear since we have
\begin{eqnarray*}
(\tilde{f}\circ\mu_M)\left[m\otimes (\sum_{k\in\mathbb{N}}a_kT^k)\right]&=&\tilde{f}(\sum_{k\in\mathbb{N}}ma_kT^k)
=\sum_{k\in\mathbb{N}}f(ma_k)T^k
=\sum_{k\in\mathbb{N}}f(m)a_kT^k\\
&=&\mu_N\left[f(m)\otimes (\sum_{k\in\mathbb{N}}a_kT^k)\right]\\
&=&(\mu_N\circ (f\otimes 1_{\tilde{A}[[T]]}))\left[m\otimes (\sum_{k\in\mathbb{N}}a_kT^k)\right],
\end{eqnarray*}
for any homogeneous morphism $f:M\longrightarrow N$ in $GR-A$ and all homogeneous elements $m\in M$ and $\sum_{k\in\mathbb{N}}a_kT^k\in\tilde{A}[[T]]$.

\vspace*{0.3cm}

(\ref{prop.the tilde functor.5}) Let $\Pi_A:\mathcal{D}(A)\longrightarrow\mathcal{H}(A)$ denote the homotopically
projective resolution functor.
Since each $e_iA$ is homotopically projective (see \cite[Example 7.6]{S}), we have that  $\Pi_A(e_iA)\cong e_iA$ in $\mathcal{H}(A)$. Moreover,  we have an isomorphism
$\mu_{e_iA}:e_iA\otimes_A\tilde{A}[[T]]\stackrel{\cong}{\longrightarrow}e_i\tilde{A}[[T]]=
e_iA\hat{\otimes}V$
in $Dg-\tilde{A}[[T]]$ (see the proof of Proposition 10.5 in \cite{S}). One then
 gets from \cite[Proposition 7.12]{S} that $\mu_{e_iA}:\mathbb{L}\iota^*(e_iA)\longrightarrow (e_iA)\uar_K^V$ is an isomorphism in
$\mathcal{D}(\tilde{A}[[T]])$, from which it follows that
$\mu_M:\mathbb{L}\iota^*(M)\longrightarrow M\uar_K^V$ is an isomorphism, for all $M\in\text{per}(A)=\text{thick}_{\mathcal{D}(A)}(e_iA\text{: }i\in I)$.
\end{proof}

\medskip

Note that we have a canonical isomorphism of dg algebras with enough idempotents
$\tilde{A}^{op}[[T]]\cong \tilde{A}[[T]]^{op}$. We will still denote by $\iota^*$
the dg functor $\tilde{A}[[T]]\otimes_A?:A-Dg\longrightarrow\tilde{A}[[T]]$ and by
$\mathbb{L}\iota^*$ its left derived functor $\mathcal{D}(A^{op})\longrightarrow\mathcal{D}(\tilde{A}[[T]]^{op})$.
We will  denote by $V\hat{\otimes}?:A-Dg\longrightarrow \tilde{A}[[T]]-Dg$ and
$\uar_k^V:\mathcal{D}(A^{op})\longrightarrow\mathcal{D}(\tilde{A}[[T]]^{op})$ the corresponding
dg functor and triangulated functor, respectively. We now get:

\begin{Cor} \label{cor.natural transformation}
Consider  the compositions of triangulated functors
$$\mathcal{D}(A)^{op}\stackrel{\mathbb{R}\text{Hom}_A(?,A)}{\longrightarrow}
\mathcal{D}(A^{op})\stackrel{\mathbb{L}\iota^*}{\longrightarrow}\mathcal{D}(\tilde{A}[[T]]^{op})$$
and  $$\mathcal{D}(A)^{op}\stackrel{(\mathbb{L}\iota^*)^o}{\longrightarrow}
\mathcal{D}(\tilde{A}[[T]])^{op}\stackrel{\mathbb{R}\text{Hom}_{\tilde{A}[[T]]}(?,\tilde{A}[[T]])}{\longrightarrow}
\mathcal{D}(\tilde{A}[[T]]^{op}). $$
There are  natural isomorphisms of triangulated functors
$$\eta :(\mathbb{L}\iota^*\circ\mathbb{R}\text{Hom}_A(?,A))_{| per(A)^{op}}\stackrel{\cong}{\longrightarrow}
[\mathbb{R}\text{Hom}_{\tilde{A}[[T]]}(?,\tilde{A}[[T]])\circ\mathbb{L}\iota^*]_{| per(A)^{op}}$$
and
$$\eta :\uar_K^V\circ\mathbb{R}\text{Hom}_A(?,A))_{| per(A)^{op}}\stackrel{\cong}{\longrightarrow}
[\mathbb{R}\text{Hom}_{\tilde{A}[[T]]}(?,\tilde{A}[[T]])\circ\uar_K^V]_{| per(A)^{op}}.$$
\end{Cor}

\begin{proof}
The first natural isomorphism is a direct consequence of
\cite[Proposition 10.5]{S}.
On the other hand, Proposition \ref{prop.the tilde functor} and its left-right
symmetric version gives natural isomorphisms
$(\mathbb{L}\iota^*)_{| per(A)}\cong (\uar_K^V)_{| per(A)}$ and
$(\mathbb{L}\iota^*)_{| per(A^{op})}\cong (\uar_K^V)_{| per(A^{op})}$. Using now the duality $$\mathbb{R}\text{Hom}_A(?,A):\text{per}(A)\stackrel{\cong^o}{\longrightarrow}\text{per}(A^{op}),$$
the result follows.
\end{proof}

\medskip

We now revisit and generalize some point of \cite{SZ}.
Note that the variable $T$ is not an element of $\tilde{A}[[T]]$, unless $A$ has a unit.
However, if $Q$ is a right dg $\tilde{A}[[T]]$-module and $x\in Q$ is a homogeneous element,
then the product $xT$ makes sense. Indeed since $x=\sum_{i\in I}xe_i$, with $xe_i=0$ for almost
all $i\in I$, the element $xT:=\sum_{i\in I}x(e_iT)$ is a well-defined element of $Q$ with
$deg(xT)=deg(x)$. Furthermore, if $f:Q\longrightarrow Q'$ is a morphism in $Dg-\tilde{A}[[T]]$,
then we have $f(xT)=f(\sum_{i\in I}x(e_iT))=\sum_{i\in I}f(x)e_iT=f(x)T$. We can now prove:

\begin{Lemma} \label{lem.t as dg transformation}
For each right dg $\tilde{A}[[T]]$-module $Q$, the map $t_Q:Q\longrightarrow Q$ ($x\rightsquigarrow xT$)
is a morphism of zero degree in $Dg-\tilde{A}[[T]]$ and, when $Q$ varies, the $t_Q$
give a homological natural transformation of dg functors
$t:1_{Dg-\tilde{A}[[T]]}\longrightarrow 1_{Dg-\tilde{A}[[T]]}$.
\end{Lemma}

\begin{proof}
Let $x\in Q$ and $\sum_{k\in\mathbb{N}}a_kT^k\in\tilde{A}[[T]]$ be homogeneous elements.
By definition of $\tilde{A}[[T]]$ and by the fact that $Q=\bigoplus_{i\in I}Qe_i$, we have a
finite subset $F\subset I$ such that $xe_i=0$ and $a_ke_i=0$, for all $i\in I\setminus F$
and all $k\in\mathbb{N}$. It follows that $(x\sum_{k\in\mathbb{N}}a_kT^k)e_i=0$, for all
$i\in I\setminus F$. We then have

\begin{eqnarray*}
t_Q(x\sum_{k\in\mathbb{N}}a_kT^k)&=&(x\sum_{k\in\mathbb{N}}a_kT^k)T
=\sum_{i\in F}(x\sum_{k\in\mathbb{N}}a_kT^k)e_iT
=\sum_{i\in F}x(\sum_{k\in\mathbb{N}}a_ke_iT^{k+1})\\
&=&x\sum_{k\in\mathbb{N}}a_kT^{k+1}
=(xT)\sum_{k\in\mathbb{N}}a_kT^k=t_Q(x)\sum_{k\in\mathbb{N}}a_kT^k,
\end{eqnarray*}
which shows that $t_Q$ is a morphism of zero degree in $Dg-\tilde{A}[[T]]$.

If $f:Q\longrightarrow Q'$ is a homogeneous morphism in $Dg-\tilde{A}[[T]]$ then
we have $(t_{Q'}\circ f)(x)=f(x)T=f(xT)=(f\circ t_Q)(x)$, for each $x\in Q$.
This proves that, when $Q$ varies, the $t_Q$ give a natural transformation of
dg functors $t:1_{Dg-\tilde{A}[[T]]}\longrightarrow 1_{Dg-\tilde{A}[[T]]}$.
This natural transformation is homological since we have
\begin{eqnarray*}
(d_Q\circ t_Q)(x)&=&d_Q(xT)=d_Q(\sum_{i\in I}x(e_iT))
=\sum_{i\in I}d_Q(x(e_iT))\\
&=&\sum_{i\in I}(d_Q(x)e_i)T=d_Q(x)T=(t_Q\circ d_Q)(x),
\end{eqnarray*}
for each homogeneous element $x\in Q$,  due to the fact that $d_{\tilde{A}[[T]]}(e_iT)=d(e_i)T=0$ (see \cite[Remark 7.1]{S}).
\end{proof}

\medskip

Note that the associated natural transformation of triangulated functors
$t:1_{\mathcal{D}(\tilde{A}[[T]])}\longrightarrow 1_{\mathcal{D}(\tilde{A}[[T]])}$
is the one given in \cite{SZ}, after translation to the language of dg
algebras with enough idempotents.
Replacing $A$ by $A^{op}$ in Lemma~\ref{lem.t as dg transformation} and interpreting right dg modules over
$\tilde{A}[[T]]^{op}\cong\tilde{A}^{op}[[T]]$ as left dg $\tilde{A}[[T]]$-modules,
we get a natural transformation of dg functors
$t:1_{\tilde{A}[[T]]-Dg}\longrightarrow 1_{\tilde{A}[[T]]-Dg}$ which in turn gives
a natural transformation of triangulated functors
$t:1_{\mathcal{D}(\tilde{A}[[T]]^{op})}\longrightarrow 1_{\mathcal{D}(\tilde{A}[[T]]^{op})}$.
These natural transformations do not have correspondents
for dg $\tilde{A}[[T]]-\tilde{A}[[T]]-$bimodules,
because the action of $T$ by multiplication on a dg
$\tilde{A}[[T]]-\tilde{A}[[T]]-$bimodule need not be the same on the left
and on the right. We say that a dg $\tilde{A}[[T]]-\tilde{A}[[T]]-$bimodule
$X$ is \emph{T-symmetric} when $Tx=xT$, for each $x\in X$.  Note that, with a
suitable modification of the argument used in the proof of
Lemma \ref{lem.t as dg transformation}, one easily sees that if $X$ is a $T$-symmetric
dg $\tilde{A}[[T]]$-bimodule, then the assignment $x\rightsquigarrow xT=Tx$ is
a morphism of $\tilde{A}[[T]]-\tilde{A}[[T]]-$bimodules, which we also denote by $t_X$.

\begin{Prop} \label{prop.natural transformation t}
Let us consider the bi-triangulated functor $$\mathbb{R}\text{HOM}_{\tilde{A}[[T]]}(?,?):
\mathcal{D}(\tilde{A}[[T]])^{op}\otimes\mathcal{D}(\tilde{A}[[T]]\otimes\tilde{A}[[T]]^{op})
\longrightarrow\mathcal{D}(\tilde{A}[[T]]^{op})$$
(see Proposition \ref{prop.bifunctors-onesidedfunctors}) and let $Q$ be a
right dg $\tilde{A}[[T]]$-module and  $X$
be a $T$-symmetric dg $\tilde{A}[[T]]-\tilde{A}[[T]]-$bimodule.

Then $\mathbb{R}\text{HOM}_{\tilde{A}[[T]]}(t_Q^o,1_X)$ and
$\mathbb{R}\text{HOM}_{\tilde{A}[[T]]}(1_Q^o,t_X)$, considered as maps $$\mathbb{R}\text{HOM}_{\tilde{A}[[T]]}(Q,X)\longrightarrow\mathbb{R}\text{HOM}_{\tilde{A}[[T]]}(Q,X)$$
are equal.
Moreover they are equal to the evaluation of the natural transformation
$$t:1_{\mathcal{D}(\tilde{A}[[T]]^{op})}\longrightarrow 1_{\mathcal{D}(\tilde{A}[[T]]^{op})}$$
at $\mathbb{R}\text{HOM}_{\tilde{A}[[T]]}(Q,X)$.
\end{Prop}

\begin{proof}
By Proposition \ref{prop.bifunctors-onesidedfunctors}, we have a natural isomorphism
of triangulated functor
$$\mathbb{R}\text{HOM}_{\tilde{A}[[T]]}(?,X)\cong\mathbb{R}\text{Hom}_{\tilde{A}[[T]]}(?,X).$$
So in order to see that $\mathbb{R}\text{HOM}_{\tilde{A}[[T]]}(t_Q^o,1_X)$ is the evaluation of $t$ at $\mathbb{R}\text{HOM}_{\tilde{A}[[T]]}(Q,X)$ it is enough to check that
$t_Q^*=\mathbb{R}\text{Hom}_{\tilde{A}[[T]]}(?,X)(t_Q)$ is precisely
$t_{\mathbb{R}\text{Hom}_{\tilde{A}[[T]]}(Q,X)}$, where
$\mathbb{R}\text{Hom}_{\tilde{A}[[T]]}(Q,X):=\mathbb{R}\text{Hom}_{\tilde{A}[[T]]}(?,X)(Q)$ in the rest of the proof.
Note that $t_Q^*$ is  the map
$$\mathbb{R}\text{Hom}_{\tilde{A}[[T]]}(Q,X)
=\overline{HOM}_{\tilde{A}[[T]]}(\Pi (Q),X)\stackrel{\Pi (t_Q)^*}{\longrightarrow}
\overline{HOM}_{\tilde{A}[[T]]}(\Pi (Q),X)=\mathbb{R}\text{Hom}_{\tilde{A}[[T]]}(Q,X).$$
Here and in the rest of the proof
$$\Pi :=\Pi_{\tilde{A}[[T]]}:\mathcal{D}(\tilde{A}[[T]])\longrightarrow\mathcal{H}(\tilde{A}[[T]])$$
and
$$\Upsilon:=\Upsilon_{\tilde{A}[[T]]\otimes\tilde{A}[[T]]^{op}}:
\mathcal{D}(\tilde{A}[[T]]\otimes\tilde{A}[[T]]^{op})\longrightarrow
\mathcal{H}(\tilde{A}[[T]]\otimes\tilde{A}[[T]]^{op})$$
are the homotopically projective and the homotopically injective resolution functors, respectively.

It is convenient to have a careful look at a special case of the action of $\Pi$ and
$\Upsilon$ on morphisms. Let $Q$ and $X$ be as in the statement and let $f:Q\longrightarrow Q$
and $\alpha :X\longrightarrow X$ be morphisms in $\mathcal{H}(\tilde{A}[[T]])$ and $\mathcal{H}(\tilde{A}[[T]]\otimes\tilde{A}[[T]]^{op})$, respectively. Abusing notation,
we put $q(f)=f$ and $q(\alpha )=\alpha$, where $q$ is the functor from the homotopy to the
derived category in each case. Viewing $Q$ and $X$ as objects of the respective derived categories,
we have a counit map $$\pi_Q:(\Pi\circ q)(Q)=\Pi (Q)\longrightarrow Q$$ and a unit map
$$\iota_X:X\longrightarrow (\Upsilon\circ q)(X)=\Upsilon (X),$$ which are quasi-isomorphism.
Then $\Pi(f):=(\Pi\circ q)(f)$ is a morphism $\Pi (Q)\longrightarrow\Pi (Q)$ in
$\mathcal{H}(\tilde{A}[[T]])$ such that $$\pi_Q\circ\Pi(f)=f\circ\pi_Q\;\; (*),$$
due to the naturality of the counit $\pi$. But since we have an isomorphism
$$\text{Hom}_{\mathcal{H}(\tilde{A}[[T]])}(\Pi (Q),\Pi(Q))
\stackrel{\cong}{\longrightarrow}\text{Hom}_{\mathcal{D}(\tilde{A}[[T]])}(Q,Q)$$
(which maps $\varphi\rightsquigarrow q(\pi )\circ \varphi \circ q(\pi )^{-1}$),
we see that $\Pi (f)$ is the unique morphism in $\mathcal{H}(\tilde{A}[[T]])$
satisfying the equality (*). Similarly, $\Upsilon (\alpha )$ is the unique morphism
$\Upsilon (X)\longrightarrow\Upsilon (X)$ in $\mathcal{H}(\tilde{A}[[T]])\otimes\tilde{A}[[T]]^{op})$
such that $$\Upsilon (\alpha )\circ\iota_X=\iota_X\circ\alpha.$$
By taking $f=t_Q$
in this argument, we readily see that $\Pi (t_Q)=t_{\Pi (Q)}$ since, due to the naturality of
$t:1_{Dg-\tilde{A}[[T]]}\longrightarrow 1_{Dg-\tilde{A}[[T]]}$,
we have that $\pi\circ t_{\Pi (Q)}=t_Q\circ\pi$ in $\mathcal{H}(\tilde{A}[[T]])$.
The analogous fact does not work for $\alpha =t_X$ since we do not have a
correspondent of the natural transformation $t$ for $\tilde{A}[[T]]-\tilde{A}[[T]]$-bimodules.
In any case, these comments together with the previous paragraph show that
$\mathbb{R}\text{Hom}_{\tilde{A}[[T]]}(?,X)(t_Q)=t_{\Pi (Q)}^*$. Using the naturality of $t$,
we see that $t_{\Pi (Q)}^*$ is a morphism
$\overline{\text{HOM}}_A(\Pi (Q),X)\longrightarrow\overline{\text{HOM}}_A(\Pi (Q),X)$
of left dg $\tilde{A}[[T]]$-modules which maps
$$f\rightsquigarrow (-1)^{|t_{\Pi(Q)}| |f|}f\circ t_{\Pi (Q)}=f\circ t_{\Pi (Q)}=t_X\circ f.$$
We then have that
$$[t_{\Pi (Q)}^*(f)](z)=(f\circ t_X)(z)=f(zT)=f(z)T=Tf(z)=(Tf)(z)=
[t_{\overline{\text{HOM}}_A(\Pi (Q),X)}(f)](z), $$
for all homogeneous elements $z\in\Pi (Q)$, using the definition of the left
$\tilde{A}[[T]]$-module structure on $\overline{\text{HOM}}_A(\Pi (Q),X)$
(see \cite[Section 8]{S}) and the $T$-symmetry of $X$. Therefore we have
$\mathbb{R}\text{Hom}_A(?,X)(t_Q)=t_{\mathbb{R}\text{Hom}_A(Q,X)}$, as desired.

On the other hand, by Proposition \ref{prop.bifunctors-onesidedfunctors}, there
is a natural isomorphism
$$\mathbb{R}\text{HOM}_{\tilde{A}[[T]]}(Q,?)\cong\mathbb{R}\text{Hom}_{\tilde{A}[[T]]}(\Pi (Q),?)=q\circ\overline{HOM}_{\tilde{A}[[T]]}(\Pi (Q),?)\circ\Upsilon$$
of triangulated functors
$$\mathcal{D}(\tilde{A}[[T]]\otimes\tilde{A}[[T]]^{op})\longrightarrow\mathcal{D}(\tilde{A}[[T]]^{op}).$$
Then $\mathbb{R}\text{HOM}_{\tilde{A}[[T]]}(1_Q^o,t_X)$ is the morphism
\begin{eqnarray*}
\overline{HOM}_{\tilde{A}[[X]]}(\Pi (Q),\Upsilon (X))&\longrightarrow&
\overline{HOM}_{\tilde{A}[[X]]}(\Pi (Q),\Upsilon (X))\\
f&\rightsquigarrow& (-1)^{|1_Q| |q(\Upsilon (t_X)|}\Upsilon (t_X)\circ f\circ 1_Q)=\Upsilon (t_X)_*(f).
\end{eqnarray*}
In other words, we have that
$$\mathbb{R}\text{HOM}_{\tilde{A}[[T]]}(1_Q^o,t_X)=q(\Upsilon (t_X)_*),$$
where $q:\mathcal{H}(\tilde{A}[[T]]^{op})\longrightarrow\mathcal{D}(\tilde{A}[[T]]^{op})$
is the canonical functor and
$$\Upsilon (t_X)_*=\overline{HOM}_A(\Pi (M),?)(\Upsilon (t_X)):\overline{HOM}_A(\Pi (M),\Upsilon (X))\longrightarrow\overline{HOM}_A(\Pi (M),\Upsilon (X)).$$
But the induced functor
$$\overline{HOM}_{\tilde{A}[[T]]}(\Pi(M),?):\mathcal{H}(\tilde{A}[[T]]\otimes\tilde{A}[[T]]^{op})
\longrightarrow
\mathcal{H}(\tilde{A}[[T]]^{op})$$
preserves quasi-isomorphisms since $\Pi (M)$ is homotopically projective in
$\mathcal{H}(\tilde{A}[[T]])$.
If now $\iota :=\iota_X:X\longrightarrow\Upsilon (X)$ is as above,  then
$$\iota_*:=\overline{HOM}_{\tilde{A}[[T]]}(\Pi (M),\iota):
\overline{HOM}_{\tilde{A}[[X]]}(\Pi (Q),X)\longrightarrow
\overline{HOM}_{\tilde{A}[[X]]}(\Pi (Q),\Upsilon (X)))$$
is a quasi-isomorphism of left dg $\tilde{A}[[T]]$-modules.
Applying to the equality $\Upsilon (t_X)\circ\iota =\iota\circ t_X$ the functor
$$\overline{HOM}_A(\Pi (M),?):\mathcal{H}(\tilde{A}[[T]]\otimes\tilde{A}[[T]]^{op})\longrightarrow
\mathcal{H}(\tilde{A}[[T]]^{op}),$$
we get the following commutative diagram in $\mathcal{H}(\tilde{A}[[T]]^{op})$,
where the horizontal arrows are quasi-isomorphisms.
$$\xymatrix{
\ol{HOM}_{{\widetilde A}[[T]]}(\Pi(M),X)\ar[r]^{\iota_*}\ar[d]_{(t_X)_*}&
\ol{HOM}_{{\widetilde A}[[T]]}(\Pi(M),\Upsilon(X))\ar[d]_{\Upsilon(t_X)_*}\\
\ol{HOM}_{{\widetilde A}[[T]]}(\Pi(M),X)\ar[r]^{\iota_*}&\ol{HOM}_{{\widetilde A}[[T]]}(\Pi(M),\Upsilon(X))
}$$
Moreover,  the left vertical arrow takes $f\rightsquigarrow t_X\circ f$,
for each homogeneous element $f\in\overline{HOM}_{\tilde{A}[[T]]}(\Pi (M),X)$.
But, in turn, we have that
$$(t_X\circ f)(v)=Tf(v)=(Tf)(v)=t_{\overline{HOM}_{\tilde{A}[[T]]}(\Pi (M),X)}(f)(v),$$
for each homogeneous element $v\in\Pi (M)$. Therefore the left vertical arrow of last
diagram is  the evaluation of the natural transformation of dg functors
$t:1_{\tilde{A}[[T]]-Dg}\longrightarrow 1_{\tilde{A}[[T]]-Dg}$ at
$\overline{HOM}_{\tilde{A}[[T]]}(\Pi (Q),X)$. The fact that $t$ is a natural transformation
of dg functors implies that we also have an equality
$$t_{\overline{HOM}_{\tilde{A}[[T]]}(\Pi (Q),\Upsilon (X))}\circ\iota_*=
\iota_*\circ t_{\overline{HOM}_{\tilde{A}[[T]]} (\Pi (Q),X)}$$
in $\tilde{A}[[T]]-Dg$
and, hence,  also in $\mathcal{H}(\tilde{A}[[T]]^{op})$. We then have that
$$t_{\overline{HOM}_{\tilde{A}[[T]]}(\Pi (Q),\Upsilon (X))}\circ\iota_*=\Upsilon (t_X)_*\circ\iota_*$$
in $\mathcal{H}(\tilde{A}[[T]]^{op})$. Applying the functor $q:\mathcal{H}(\tilde{A}[[T]]^{op})\longrightarrow\mathcal{D}(\tilde{A}[[T]]^{op})$
to this last equality and bearing in mind that $q(\iota_*)$ is an isomorphism,
we conclude that $$q((\Upsilon (t_X))_*)=q(t_{\overline{HOM}_{\tilde{A}[[T]]}(\Pi (Q),\Upsilon (X))})=t_{\mathbb{R}\text{HOM}_{\tilde{A}[[T]]}(Q,X)}.$$
\end{proof}

\medskip

For our next result we adopt the terminology of \cite[Proposition 9]{SZ} and, for
the given dg algebra with enough idempotents $A$, we put
$\mathcal{C}_V^o=\text{per}(\tilde{A}[[T]])$,
we denote by $\mathcal{C}_V^o[t^{-1}]$ the localization of $\mathcal{C}_V^o$
with respect to natural
transformation $t$ given above (see \cite[Remark 2]{SZ} for the definition) and we let
$p:\mathcal{C}_V^o\longrightarrow\mathcal{C}_V^o[t^{-1}]$ be the canonical functor.
We also put
${}_V\mathcal{C}^o=\text{per}(\tilde{A}[[T]]^{op})$, ${}_V\mathcal{C}^o[t^{-1}]$ and
$p:{}_V\mathcal{C}^o\longrightarrow {}_V\mathcal{C}^o[t^{-1}]$
for the corresponding concepts on the left.

\begin{Lemma} \label{lem.transformation t versus duality}
Let $p:\mathcal{C}_V^o\longrightarrow\mathcal{C}_V^o[t^{-1}]$ and
$p':{}_V\mathcal{C}^o\longrightarrow {}_V\mathcal{C}^o[t^{-1}]$ be the canonical
triangulated functors given by localization, and let
$Q_1$ and  $Q_2$  be objects of $\mathcal{C}_V^o=\text{per}(\tilde{A}[[T]])$.
There is an isomorphism $p(Q_1)\cong p(Q_2)$ if, and only if, there is an isomorphism
$p'(Q_1^\star )\cong p'(Q_2^\star)$,
where $(?)^\star:=\mathbb{R}\text{Hom}_{\tilde{A}[[T]]}(?,\tilde{A}[[T]]):
\mathcal{D}(\tilde{A}[[T]])^{op}=\mathcal{C}_V^{op}\longrightarrow
{}_V\mathcal{C}=\mathcal{D}(\tilde{A}^{op}[[T]])$ is the usual triangulated functor.
\end{Lemma}

\begin{proof}
The fact that $p(Q_1)$ and $p(Q_2)$ are isomorphic in $\mathcal{C}_V^o[t^{-1}]$ means
that we have morphisms $f:Q_1\longrightarrow Q_2$ and $g:Q_2\longrightarrow Q_1$ in
$\mathcal{C}_V^o=\text{per}(\tilde{A}[[T]])$ such that $g\circ f\circ t_{Q_1}^r=t_{Q_1}^s$
and $f\circ g\circ t_{Q_2}^m=t_{Q_2}^n$, for some  $r,s,m,n\in\mathbb{N}$.
If we now apply the duality $$(?)^\star=\mathbb{R}\text{Hom}_{\tilde{A}[[T]]}(?,\tilde{A}[[T]]):
\text{per}(\tilde{A}[[T]])=\mathcal{C}_V^o\stackrel{\cong^o}{\longrightarrow}{}_V\mathcal{C}^o=
\text{per}(\tilde{A}[[T]]^{op}),$$ then we get that $(t_{Q_2}^\star )^m\circ g^\star\circ f^\star
=(t_{Q_2}^\star )^n$.
But
Propositions \ref{prop.bifunctors-onesidedfunctors} and \ref{prop.natural transformation t}
tell us that we have
$(t_{Q_k})^\star =t_{Q_k^\star}$ for $k=1,2$, which implies that
$p'(g^\star)\circ p'(f^\star)$ and $p'(f^\star)\circ p'(g^\star)$ are
isomorphisms in ${}_V\mathcal{C}^o[t^{-1}]$, and hence that
$p'(Q_1^\star)\cong p'(Q_2^\star )$ in ${}_V\mathcal{C}^o[t^{-1}]$.

The reverse implication follows by exchanging the roles of $A$ and
$A^{op}$ and of $Q_k$ and $Q_k^\star$, bearing in mind that $Q_k$ is
isomorphic to $Q_k^{\star\star}$, for $k=1,2$ (see \cite[Proposition 10.4]{S}).
\end{proof}

\medskip

The first assertion of the following Lemma~\ref{lem.perfects image of thick} seems to be folklore,
but we include a
short proof. A right dg $\tilde{A}[[T]]$-module $Q$ will be called
\emph{T-torsion-free} when $t_Q:Q\longrightarrow Q$ is monomorphism in
$Gr-\tilde{A}[[T]]$. Note that this is equivalent to saying that $t_Q$ is a
monomorphism for the abelian structure of $\mathcal{C}(\tilde{A}[[T]])$.

\begin{Lemma} \label{lem.perfects image of thick}
Let $A$ be a dg algebra with enough idempotents and let
$q_A:\mathcal{H}(A)\longrightarrow\mathcal{D}(A)$  be the canonical functor. The following assertions hold:
\begin{enumerate}
\item The induced functor
$q:\text{thick}_{\mathcal{H}(A)}(e_iA\text{: }i\in I)\longrightarrow\text{per}(A)=\mathcal{D}^c(A)$
is an equivalence of triangulated categories.
\item Each $Q\in\text{thick}_{\mathcal{H}(\tilde{A}[[T]])}(e_i\tilde{A}[[T]]\text{: }i\in I)$
is isomorphic in $\mathcal{H}(\tilde{A}[[T]])$ to a T-torsion-free  right dg $\tilde{A}[[T]]$-module.
\end{enumerate}
\end{Lemma}

\begin{proof}
(1) The subcategory $\text{thick}_{\mathcal{H}(A)}(e_iA\text{: }i\in I)$ of $\mathcal{H}(A)$
consists of homotopically projective objects and the restriction of $q$ to the subcategory of
homotopically projective objects is fully faithful. In order to prove the density, recall that
$\text{per}(A)=\text{thick}_{\mathcal{D}(A)}(e_iA\text{: }i\in I)$ (see \cite[Theorem 5.3]{K1}).
This implies in particular that each  $X\in\text{per}(A)$ is a direct summand in $\mathcal{D}(A)$
of a right dg $A$-module $P$ for which there is a sequence of morphisms $$0=P_0\stackrel{f_1}{\longrightarrow}P_1\stackrel{f_2}{\longrightarrow}\cdots\stackrel{f_{n-1}}{\longrightarrow}
P_{n-1}\stackrel{f_n}{\longrightarrow}P_n$$
in $\mathcal{D}(A)$ such that $P_n=P$ and  $\text{cone}(f_k)\cong e_{i_k}A[m_k]$,
for some  $i_k\in I$ and some $m_k\in\mathbb{Z}$, for $k=1,...,n$.
We will prove by induction on $n$ that $P\cong q(Q)$,
for some $Q\in\text{thick}_{\mathcal{H}(A)}(e_iA\text{: }i\in I)$. For $n=0$ there is
nothing to prove, so we assume that $n>1$. By the induction hypothesis, we can choose
$Q_{n-1}\in\text{thick}_{\mathcal{H}(A)}(e_iA\text{: }i\in I)$ such that $q(Q_{n-1})\cong P_{n-1}$.
We then get a distinguished triangle
$Q_{n-1}\longrightarrow P\longrightarrow e_iA[m]\stackrel{f[1]}{\longrightarrow}Q_{n-1}[1]$,
for some $i\in I$, some $m\in\mathbb{Z}$ and some morphism
$f:e_iA[m]\longrightarrow Q_{n-1}$ in $\mathcal{D}(A)$. But the functor $q$ gives an isomorphism $\text{Hom}_{\mathcal{H}(A)}(e_iA,Q_{n-1})\stackrel{\cong}{\longrightarrow}
\text{Hom}_{\mathcal{D}(A)}(e_iA,Q_{n-1})$.
This means that we may view $f$ as a morphism in $\mathcal{H}(A)$, and then the triangulated  cone
$Q=\text{cone}_{\mathcal{H}(A)}(f)$ is in $\text{thick}_{\mathcal{H}(A)}(e_iA\text{: }i\in I)$
and satisfies that $q(Q)\cong P$.

Let  now $X$, $P$ and $Q$ be as above and let $e\in\text{End}_{\mathcal{D}(A)}(P)$
be the idempotent endomorphism corresponding to the direct summand $X$ of $P$. Since $q$
gives an algebra isomorphism
$\text{End}_{\mathcal{H}(A)}(Q)\stackrel{\cong}{\longrightarrow}\text{End}_{\mathcal{D}(A)}(P)$,
we have a unique $\epsilon =\epsilon^2\in\text{End}_{\mathcal{H}(A)}(Q)$ such that $q(\epsilon )=e$.
Since $\mathcal{H}(A)$ has arbitrary (set-indexed) coproducts, we know that idempotents split
in $\mathcal{H}(A)$ (see \cite[Proposition 1.6.8]{Neeman}). We then get a direct summand $Y$ of $Q$
in $\mathcal{H}(A)$ corresponding to $\epsilon$, and we clearly have that $q(Y)\cong X$.

\vspace*{0.3cm}

(2) Let $Q\in\text{thick}_{\mathcal{H}(\tilde{A}[[T]])}(e_i\tilde{A}[[T]]\text{: }i\in I)$
be any object. By the obvious adaptation of \cite[Theorem 3.1]{K1} to the language of
dg algebras with enough idempotents, we know that there is a chain of inflations in
$\mathcal{C}(\tilde{A}[[T]])$  $$0=P_0\hookrightarrow P_1\hookrightarrow ...
\hookrightarrow P_n\hookrightarrow ... $$
such that  $\text{Coker}(P_{n-1}\hookrightarrow P_n)$ is a direct summand in
$\mathcal{C}(\tilde{A}[[T]])$ of a (possibly infinite) coproduct of dg right
$\tilde{A}[[T]]$-module of the form $e_i\tilde{A}[[T]] [m]$, with $i\in I$ and
$m\in\mathbb{Z}$, and $P=\bigcup_{n\in\mathbb{N}}P_n$ is isomorphic to $Q$ in
$\mathcal{D}(\tilde{A}[[T]])$. Since all the exact sequences
$$0\rightarrow P_{n-1}\hookrightarrow P_n\longrightarrow P_n/P_{n-1}\rightarrow 0$$
split in $Gr-\tilde{A}[[T]]$, we readily see that $P$ is projective in this category.
In particular $P$ is $T$-torsion-free. But $P$ and $Q$ are homotopically
projective objects of $\mathcal{H}(\tilde{A}[[T]])$, which implies that the canonical functor
$q:\mathcal{H}(\tilde{A}[[T]])\longrightarrow\mathcal{D}(\tilde{A}[[T]])$ induces bijections
$\text{Hom}_{\mathcal{H}(\tilde{A}[[T]])}(X,Y)\stackrel{\cong}{\longrightarrow}
\text{Hom}_{\mathcal{D}(\tilde{A}[[T]])}(X,Y)$,
for $X,Y\in\{P,Q\}$. We deduce that any isomorphism
$P\stackrel{\cong}{\longrightarrow}Q$ in $\mathcal{D}(\tilde{A}[[T]])$ can be
lifted to a corresponding isomorphism in $\mathcal{H}(\tilde{A}[[T]])$.
\end{proof}

\subsection{The main theorem under hypothesis (b)}

We can now complete the proof of Theorem \ref{maintheoremintroduction}.

\begin{Theorem} \label{thm.symmetry of degenerations}
Let $\mathcal{C}_k^0$ be the category of compact objects of an algebraic
compactly generated triangulated category. For any objects
$M,N\in\text{Ob}(\mathcal{C}_k^0)$, the following assertions are equivalent:
\begin{enumerate}
\item There is a distinguished triangle
$Z_\ell\stackrel{\begin{pmatrix} v\\ u\end{pmatrix}}{\longrightarrow}Z_\ell\oplus M
\stackrel{\begin{pmatrix} h & j \end{pmatrix}}{\longrightarrow} N\longrightarrow Z_\ell[1]$,
where $v$ is a nilpotent endomorphism of $Z_\ell$.
\item There is a distinguished triangle
$N\stackrel{\begin{pmatrix} j\\ h \end{pmatrix}}{\longrightarrow} M\oplus Z_r
\stackrel{\begin{pmatrix} u & v\end{pmatrix}}{\longrightarrow} Z_r\longrightarrow N[1]$,
where $v$ is a nilpotent endomorphism of $Z_r$.
\end{enumerate}
\end{Theorem}

\begin{proof}
Using the version of Keller's theorem for dg algebras with
enough idempotents (see \cite[Corollary 6.11]{S}), we can and shall assume that
$\mathcal{C}_k^o=\mathcal{D}(A)^c=\text{per}(A)$, where $A$ is a dg algebra
with enough idempotents.

\vspace*{0.3cm}

$(1)\Longrightarrow (2):$ In  \cite[Proposition 9]{SZ} we showed that if there is a
distinguished triangle as in assertion 1, then the quintuple
$(\mathcal{C}_k,\mathcal{C}_V,\mathcal{C}_V^o,\uar_k^V,t)$ give degeneration
data for $\mathcal{C}_k^o$, where
$\mathcal{C}_k=\mathcal{D}(A)$, $\mathcal{C}_V=\mathcal{D}(\tilde{A}[[T]])$, $\mathcal{C}_V^o=\mathcal{D(\tilde{A}[[T]])}^c=\text{per}(\tilde{A}[[T]])$ and
$\uar_k^V:\mathcal{D}(A)\longrightarrow\mathcal{D}(\tilde{A}[[T]])$ and
$t:1_{\mathcal{D}(\tilde{A}[[T]])}\longrightarrow 1_{\mathcal{D}(\tilde{A}[[T]])}$
are as in the previous results of this section. Moreover, in the above mentioned result
\cite[Proposition 9]{SZ} it was also proved that there exists an object
$Q\in\mathcal{C}_V^o=\text{per}(\tilde{A}[[T]])$ so that both required conditions
for categorical degeneration are satisfied, namely:

\begin{enumerate}
\item If $p:\mathcal{C}_V^o\longrightarrow\mathcal{C}_V^o[t^{-1}]$ is the
canonical functor, then $p(Q)\cong p(M\uar_k^V)$;
\item $\phi (\text{cone}(t_Q))\cong N$, where
$\phi :\mathcal{C}_V^o=\text{per}(\tilde{A}[[T]])\longrightarrow\mathcal{D}(A)=\mathcal{C}_k$
is the restriction  of $\iota_*:\mathcal{D}(\tilde{A}[[T]])\longrightarrow\mathcal{D}(A)$
to $\text{per}(\tilde{A}[[T]])$.
\end{enumerate}
Here and in the rest of the proof $\text{cone}(f)$ denotes the
triangulated cone. With this information in mind, we give the proof of the
theorem, which is divided in two steps:

\emph{Step 1: If $Q_1$ is a T-torsion-free right dg $\tilde{A}[[T]]$-module in
$\text{thick}_{\mathcal{H}(\tilde{A}[[T]])}(e_i\tilde{A}[[T]]\text{: }i\in I)$ and if we put
$Q_1^\star =\overline{HOM}_{\tilde{A}[[T]]}(Q_1,\tilde{A}[[T]])$, then
$$\phi (\text{cone}(t_{Q_1^\star}))\cong \mathbb{R}\text{Hom}_A(?,A)(\phi (\text{cone}(t_{Q_1})))$$}

Note that $Q_1$ is homotopically projective, so that we also have
$$Q_1^\star
=\mathbb{R}\text{Hom}_{\tilde{A}[[T]]}(?,\tilde{A}[[T]])(Q_1)\cong
\mathbb{R}\text{HOM}_{\tilde{A}[[T]]}(Q_1,\tilde{A}[[T]])$$
(see Proposition \ref{prop.bifunctors-onesidedfunctors}).
On the other hand, the homomorphism of dg algebras
$$\rho\otimes\rho^o :\tilde{A}[[T]]\otimes\tilde{A}[[T]]^{op}\longrightarrow A\otimes A^{op}$$
gives a restriction of scalars functor
$$(\rho\otimes\rho^o)_*:A-Dg-A\longrightarrow\tilde{A}[[T]]-Dg-\tilde{A}[[T]].$$
In particular $A$ is a dg $\tilde{A}[[T]]-\tilde{A}[[T]]-$bimodule by defining
$(\sum_{k\in\mathbb{N}}a_kT^k)a=a_0a$ and $a (\sum_{k\in\mathbb{N}}a_kT^k))=aa_0$,
for all homogeneous elements $\sum_{k\in\mathbb{N}}a_kT^k\in\tilde{A}[[T]]$ and $a\in A$.
Note that we then have an exact sequence of T-symmetric $\tilde{A}[[T]]-\tilde{A}[[T]]-$bimodules
$$0\rightarrow\tilde{A}[[T]]\stackrel{t_{\tilde{A}[[T]]}}{\longrightarrow}\tilde{A}[[T]]
\stackrel{\rho}{\longrightarrow}A\rightarrow 0$$ in $Gr-(\tilde{A}[[T]]\otimes\tilde{A}[[T]]^{op})$
and in (the abelian structure of) $\mathcal{C}(\tilde{A}[[T]]\otimes\tilde{A}[[T]]^{op})$.

The last sequence gives a distinguished triangle $$\tilde{A}[[T]]\stackrel{t_{\tilde{A}[[T]]}}{\longrightarrow}\tilde{A}[[T]]\stackrel{\rho}{\longrightarrow}
A\longrightarrow\tilde{A}[[T]] [1]$$
in $\mathcal{D}(\tilde{A}[[T]]\otimes\tilde{A}[[T]]^{op})$.
By   Propositions \ref{prop.bifunctors-onesidedfunctors} and \ref{prop.natural transformation t},
application of the functor
$\mathbb{R}\text{Hom}_{\tilde{A}[[T]]}(Q_1,?):\mathcal{D}(\tilde{A}[[T]]\otimes\tilde{A}[[T]]^{op})
\longrightarrow
\mathcal{D}(\tilde{A}[[T]]^{op})$
to the last distinguished triangle gives a distinguished triangle
$$Q_1^\star\stackrel{t_{Q_1^\star}}{\longrightarrow}Q_1^\star\longrightarrow
\mathbb{R}\text{Hom}_{\tilde{A}[[T]]}(Q_1,A)\longrightarrow Q_1^\star [1]$$ in
$\mathcal{D}(\tilde{A}[[T]]^{op})$, so that
$\text{cone}(t_{Q_1^\star})\cong\mathbb{R}\text{Hom}_{\tilde{A}[[T]]}(Q_1,A):=
\mathbb{R}\text{Hom}_{\tilde{A}[[T]]}(Q_1,?)(A)$.

It is important to notice that, by Proposition~\ref{prop.bifunctors-onesidedfunctors}
again, we have  isomorphisms in $\mathcal{D}(\tilde{A}[[T]]^{op})$:
$$\mathbb{R}\text{Hom}_{\tilde{A}[[T]]}(Q_1,?)(A)
\cong\mathbb{R}\text{HOM}_{\tilde{A}[[T]]}(Q_1,A)\cong
\mathbb{R}\text{Hom}_{\tilde{A}[[T]]}(?,A)(Q_1). $$
When we apply the contravariant triangulated functor $\mathbb{R}\text{Hom}_{\tilde{A}[[T]]}(?,A):\mathcal{D}(\tilde{A}[[T]])\longrightarrow
\mathcal{D}(\tilde{A}[[T]]^{op})$ to the morphisms $t_{Q_1}:Q_1\longrightarrow Q_1$
we obtain the zero map. Indeed, due to the homotopically projective condition of $Q_1$,
we have that  $\mathbb{R}\text{Hom}_{\tilde{A}[[T]]}(?,A)(Q_1)=
\overline{HOM}_{\tilde{A}[[T]]}(Q_1,A)$.
But since multiplication by $T$ kills the elements of $A$, for each homogeneous element
$f\in\overline{HOM}_{\tilde{A}[[T]]}(Q_1,A)$ we have
$$[\mathbb{R}\text{Hom}_{\tilde{A}[[T]]}(?,A)(t_{Q_1})](f)=
(-1)^{|f| |t_{Q_1}|}f\circ t_{Q_1}=f\circ t_{Q_1},$$
and this is a morphism of right dg $\tilde{A}[[T]]$-modules  $Q_1\longrightarrow A$
such that $$(f\circ t_{Q_1})(x)= f(xT)=f(x)T=0$$ for all $x\in Q_1.$ On the other hand,
the contravariant dg functor
$\overline{HOM}_{\tilde{A}[[T]]}(?,A):Gr-\tilde{A}[[T]]\longrightarrow
\tilde{A}[[T]]-Gr$ is left exact and we have an exact sequence
$$0\rightarrow Q_1\stackrel{t_{Q_1}}{\longrightarrow}Q_1
\stackrel{p}{\longrightarrow} Q_1/TQ_1\rightarrow 0$$
in $Gr-\tilde{A}[[T]]$ (which is actually an exact sequence in $\mathcal{C}(\tilde{A}[[T]])$)
due to the T-torsion-free condition of $Q_1$. It follows that we have an exact sequence
$$0\rightarrow\overline{HOM}_{\tilde{A}[[T]]}(Q_1/TQ_1,A)\stackrel{p^*}{\longrightarrow}
\overline{HOM}_{\tilde{A}[[T]]}(Q_1,A)\stackrel{0}{\longrightarrow}
\overline{HOM}_{\tilde{A}[[T]]}(Q_1,A)$$
in $\tilde{A}[[T]]-Gr$. Therefore we have an isomorphism
$$p^*:\overline{HOM}_{\tilde{A}[[T]]}(Q_1/TQ_1,A)
\stackrel{\cong}{\longrightarrow}
\overline{HOM}_{\tilde{A}[[T]]}(Q_1,A).$$
We will see that $p^*$
commutes with the differentials, which will show that we have an
isomorphism
$$\text{cone}(t_{Q_1^\star})\cong\overline{HOM}_{\tilde{A}[[T]]}(Q_1/TQ_1,A)$$ in
$\mathcal{D}(\tilde{A}[[T]]^{op})$. Indeed if
$d:\overline{HOM}_{\mathcal{A}[[T]]}(Q_1/Q_1T,A)\longrightarrow
\overline{HOM}_{\mathcal{A}[[T]]}(Q_1/Q_1T,A)$
and $\delta :\overline{HOM}_{\mathcal{A}[[T]]}(Q_1,A)\longrightarrow
\overline{HOM}_{\mathcal{A}[[T]]}(Q_1,A)$
denote the respective differentials, then, for each homogeneous element
$g\in\overline{HOM}_{\mathcal{A}[[T]]}(Q_1/Q_1T,A)$, we have
$$(\delta\circ p^*)(g)=\delta (g\circ p)=d_A\circ g\circ p-(-1)^{|g\circ p|}g\circ p\circ d_{Q_1}
=d_A\circ g\circ p-(-1)^{|g|}g\circ p\circ d_{Q_1}.  $$
But $p$ is a morphism in
$Z^0(\tilde{A}[[T]]-Dg)=\mathcal{C}(\tilde{A}[[T]]^{op})$, so that
$d_{Q_1/Q_1T}\circ p-p\circ d_{Q_1}=0$. We then get
$$(\delta\circ p^*)(g)=d_A\circ g\circ p-(-1)^{|g|}g\circ d_{Q_1/Q_1T}\circ p
=d(g)\circ p=(p^*\circ d)(g), $$
and hence $\delta\circ p^*=p^*\circ d$ as desired.

If we apply the restriction of scalars functor
$\phi=\iota_*:\tilde{A}[[T]]-Dg\longrightarrow A-Dg$,
then $\overline{HOM}_{\tilde{A}[[T]]}(Q_1/TQ_1,A)$ is taken to
$\overline{HOM}_A(\phi (Q_1/Q_1T),A)$.
But, as right dg $A$-modules,  we have an isomorphism
$\phi (Q_1/Q_1T)=\rho^*(Q_1)$ where
$\rho^*=?\otimes_{\tilde{A}[[T]]}A:Dg-\tilde{A}[[T]]\longrightarrow Dg-A$
is the extension of scalars
along the morphism of dg algebras with enough idempotents
$\rho :\tilde{A}[[T]]\longrightarrow A$.
Then the induced triangulated functor
$$\rho^*=H^0(\rho^*):H^0(Dg-\tilde{A}[[T]])=\mathcal{H}(\tilde{A}[[T]])
\longrightarrow\mathcal{H}(A)=H^0(Dg-A)$$
has the property that
$$\rho^*(Q_1)=\phi (Q_1/Q_1T)\in\text{thick}_{\mathcal{H}(A)}(\rho^*(e_i\tilde{A}[[T]])
\text{: }i\in I)=\text{thick}_{\mathcal{H}(A)}(e_iA\text{: }i\in I).$$
In particular, we get that $\phi (Q_1/Q_1T)$ is homotopically projective in
$\mathcal{H}(A)$ and, hence, that
$$\overline{HOM}_A(\phi (Q_1/Q_1T),A)\cong
\mathbb{R}\text{Hom}_A(?,A)(\phi (Q_1/Q_1T))\cong\mathbb{R}\text{HOM}_A(\phi (Q_1/Q_1T),A).$$
(see Proposition \ref{prop.bifunctors-onesidedfunctors}).
Bearing in mind that the exact sequence
$$0\rightarrow Q_1\stackrel{t_{Q_1}}{\longrightarrow}Q_1\longrightarrow Q_1/TQ_1
\rightarrow 0$$
in $\mathcal{C}(\tilde{A}[[T]])$
(with respect to the abelian exact structure of $\mathcal{C}(\tilde{A}[[T]])$)
gives a distinguished triangle
$$Q_1\stackrel{t_{Q_1}}{\longrightarrow}Q_1\longrightarrow Q_1/TQ_1\longrightarrow Q_1[1]$$
in $\mathcal{D}(\tilde{A}[[T]])$, we get that
$Q_1/Q_1T\cong\text{cone}(t_{Q_1})$ in $\mathcal{D}(\tilde{A}[[T]])$.
We then get isomorphisms
\begin{eqnarray*}
\phi (\text{cone}(t_{Q_1^\star}))&\cong&\phi (\overline{HOM}_{\tilde{A}[[T]]}(Q_1/TQ_1,A))\\
&\cong&
\overline{HOM}_A(\phi (Q_1/Q_1T),A)\\
&\cong& \mathbb{R}\text{Hom}_A(?,A)(\phi (\text{cone}(t_{Q_1}))).
\end{eqnarray*}

\vspace*{0.3cm}

\emph{Step 2: End of the proof:}
Let now $M$ and $N$ be as in assertion (1)
and let $Q\in\mathcal{C}_V^o=\text{per}(\tilde{A}[[T]])$ be the
right dg $\tilde{A}[[T]]$-module considered in the second paragraph
of this proof. Using Lemma \ref{lem.perfects image of thick}, without loss of
generality, we may assume that $M,N\in\text{thick}_{\mathcal{H}(A)}(e_iA\text{: }i\in I)$ and
that $Q$ is a T-torsion-free right
$\tilde{A}[[T]]$-module in $\text{thick}_{\mathcal{H}(\tilde{A}[[T]])}(e_i\tilde{A}[[T]]\text{: }i\in I)$.
Note that then
$\tilde{M}[[T]]\in\text{thick}_{\mathcal{H}(\tilde{A}[[T]])}(e_i\tilde{A}[[T]]\text{: }i\in I)$
and $\tilde{M}[[T]]$ is T-torsion-free.

According to the Step 1, we  have  that
$$\phi (\text{cone}(t_{Q^\star}))\cong\mathbb{R}\text{Hom}_A(?,A)(\phi (\text{cone}(t_{Q})))=
\mathbb{R}\text{Hom}_A(?,A)(N)=:N^*.$$
On the other hand, by Corollary \ref{cor.natural transformation}
and Proposition \ref{prop.the tilde functor}, we know that
$$(M\uar_k^V)^\star :=\mathbb{R}\text{Hom}_{\tilde{A}[[T]]}(?,\tilde{A}[[T]])(M\uar_k^V)\cong M^*\uar_k^V,$$
where $M^*:=\mathbb{R}\text{Hom}_A(?,A)(M)$. Moreover, by Lemma \ref{lem.transformation t versus duality},
we get that $p' ((M\uar_k^V)^\star)$ and $p'(Q^\star )$ are isomorphic in ${}_V\mathcal{C}^o[t^{-1}]$,
where $p':{}_V\mathcal{C}^o\longrightarrow{}_V\mathcal{C}^o[t^{-1}]$ is the canonical functor.
It follows from this that the left dg $\tilde{A}[[T]]$-module
$$Q^\star:=\mathbb{R}\text{Hom}_{\tilde{A}[[T]]}(?,\tilde{A}[[T]])(Q)=
\mathbb{R}\text{HOM}_{\tilde{A}[[T]]}(Q,\tilde{A}[[T]])$$
defines a categorical degeneration $M^*\leq_{cdeg}N^*$, where
$(?)^*=\mathbb{R}\text{Hom}_A(?,A):\text{per}(A)\stackrel{\cong^o}{\longrightarrow}
\text{per}(A^{op})$ is the duality defined by the regular dg bimodule $X=A$
(see \cite[Proposition 10.4]{S}).

But categorical degeneration implies triangle degeneration by Theorem~\ref{mainSZ} (see \cite[Proposition 8]{SZ}),
so that we have $M^*\leq_{\Delta+\text{nil}}N^*$. That is,  we have a distinguished triangle
$$U\stackrel{\begin{pmatrix} v'\\ u' \end{pmatrix}}{\longrightarrow}
U\oplus M^*\stackrel{\begin{pmatrix} h' & j' \end{pmatrix}}{\longrightarrow}N^*\longrightarrow U[1]$$
in ${}_k\mathcal{C}^o=\text{per}(A^{op})$, where $v'$ is a nilpotent endomorphism of $U$.
Applying the duality $\mathbb{R}\text{Hom}_{A^{op}}(?,A):\text{per}(A^{op})\stackrel{\cong^o}{\longrightarrow}
\text{per}(A)$ to this last distinguished triangle, we obtain a distinguished triangle
$$N\stackrel{\begin{pmatrix} j'* \\ h'^*\end{pmatrix}}{\longrightarrow}
M\oplus U^*\stackrel{\begin{pmatrix} u'^* & v'^*\end{pmatrix}}{\longrightarrow}
U^*\longrightarrow N[1],$$
where $v'^*$ is clearly a
nilpotent endomorphism of $U^*$. The proof of the implication ends by taking $Z_r:=U^*$.

\vspace*{0.3cm}

$(2)\Longrightarrow (1):$ By applying the duality
$\mathbb{R}\text{Hom}_{A}(?,A):\mathcal{C}_k^o=\text{per}(A)
\stackrel{\cong^o}{\longrightarrow}
\text{per}(A^{op})$ to the distinguished triangle in assertion (2)
we get that $M^*\leq_{\Delta+\text{nil}}N^*$ in
${}_k\mathcal{C}^o:=\text{per}(A^{op})$. The proof of the
implication $1)\Longrightarrow 2)$,
when applied with $A^{op}$ instead of $A$, shows
that we have $M\cong M^{**}\leq_{\Delta+\text{nil}}N^{**}\cong N$,
so that the distinguished triangle of assertion 1) exists.
\end{proof}

Observe that the proof actually gives the following additional statement.

\begin{Cor}
Let $\mathcal{C}_k^0={\mathcal D}(A)^c$ be the category of compact objects in the derived
category ${\mathcal D}(A)$ of a dg algebra $A$ with enough idempotents. For any objects
$M,N\in\text{Ob}(\mathcal{C}_k^0)$, we get $$M\leq_{\Delta+\text{nil}}N\Leftrightarrow
\mathbb{R}\text{Hom}_{A}(?,A)(M)\leq_{\Delta+\text{nil}}\mathbb{R}\text{Hom}_{A}(?,A)(N).$$
\end{Cor}


%
%



\begin{thebibliography}{88}

\bibitem{B}
Theo B\"uhler, {\em Exact categories}, Expositiones Mathematicae {\bf 28}(1)  (2010) 1-69

\bibitem{CYZ}
Xiao-Wu Chen, Yu Ye and Pu Zhang, {\em Algebras of derived dimension zero},
Communications in Algebra {\bf 36} (2008) 1-10.

\bibitem{H}
Dieter Happel, {\em Triangulated Categories in the Representation Theory
of Finite Dimensional Algebras}, London Mathematical Society Lecture Note Series {\bf 119}.
Cambridge University Press 1988.

\bibitem{JSZ} Bernt Tore Jensen, Xiuping Su and Alexander Zimmermann,
{\em Degeneration-like orders for triangulated categories}, Journal of Algebra
and Applications {\bf 4} (2005) 587-597.

\bibitem{K1} Bernhard Keller, {\em Deriving DG categories}, Annales
Scientifiques de l'\'Ecole Normale Sup\'erieure {\bf 27} (1994) 63-102.

\bibitem{K2} Bernhard Keller, {\em On differential graded categories}, In:
International Congress of Mathematics, vol. II. European Mathematical Society Zurich (2006) 151-190.

\bibitem{KellerScherotzke}
Bernhard Keller and Sarah Scherotzke, {\em Graded quiver varieties and derived categories},
Journal f\"ur die reine und angewandte Mathematik {\bf 713} (2016) 85-127.

\bibitem{Neeman}
Amnon Neeman, {\em Triangulated Categories}, Princeton University Press 2001.


\bibitem{Riedtmann1}
Christine Riedtmann, {\em Degenerations for representations of quivers with relations},
Annales Scientifiques de l'\'Ecole Normale Sup\'erieure {\bf 19} (1986) 275-301.

\bibitem{S}
Manuel Saor\'in, {\em Dg algebras with enough idempotents, their dg
modules and their derived categories}, Preprint available at https://arxiv.org/abs/1612.04719 (2016).

\bibitem{SZ}
Manuel Saor\'\i n and Alexander Zimmermann, {\em An Axiomatic Approach for Degenerations
in Triangulated Categories}, Applied Categorical Structures  {\bf 24} (2016) no 4, 385-405.

\bibitem{St}
Bo Stenstr\"om, {\sc Rings of quotients}, Springer-Verlag 1975.

\bibitem{Verdier}
Jean-Louis Verdier,  {\em Des cat\'egories d\'eriv\'ees des cat\'egories abeliennes}.
Ast\'erisque {\bf 239}, Soci\'et\'e Math\'ema\-tique de France (1996).

\bibitem{Yoshinomodules}
Yuji Yohino, {\em On degeneration of modules}, Journal of
Algebra {\bf 278} (2004) 217-226.

\bibitem{Yoshino}
Yuji Yoshino, {\em Stable degeneration of Cohen-Macaulay modules},
Journal of Algebra {\bf 332} (2011) 500-521.

\bibitem{Wang}
Zhengfang Wang, {\em Triangle order $\leq_\Delta$ in Singular Categories},
Algebras and Representation Theory {\bf 19} (2016) 397-404.

\bibitem{Webb}
Peter Webb, {\em Consequences of the existence of Auslander-Reiten triangles
with applications to perfect complexes for self-injective algebras}, preprint 2014;

\bibitem{reptheobuch}
Alexander Zimmermann, {\sc Representation Theory; A homological algebra point of view},
Springer Verlag, Cham 2014.

\bibitem{Zwara}
Grzegorz Zwara, {\em A degeneration-like order for modules},
Archiv der Mathematik {\bf 71} (1998) 437-444.

\bibitem{Zwara1}
Grzegorz Zwara, {\em Degenerations of finite dimensional modules
are given by extensions}, Compositio Mathematica {\bf 121} (2000) 205-218.


\end{thebibliography}
\end{document}